\documentclass[11pt]{article}

\usepackage[margin=1in]{geometry}
\usepackage[T1]{fontenc}
\usepackage[utf8]{inputenc}
\usepackage{lmodern}
\usepackage{microtype}

\usepackage{amsmath,amssymb,amsfonts,amsthm}
\usepackage{mathtools}
\usepackage{bm}

\theoremstyle{plain}
\newtheorem{theorem}{Theorem}[section]
\newtheorem{lemma}[theorem]{Lemma}
\newtheorem{proposition}[theorem]{Proposition}

\theoremstyle{definition}
\newtheorem{definition}[theorem]{Definition}

\theoremstyle{remark}
\newtheorem{remark}[theorem]{Remark}

\usepackage{graphicx}
\usepackage{float}
\usepackage{caption}
\usepackage{booktabs}
\usepackage{array}
\graphicspath{{figures/}}

\usepackage{xcolor}
\usepackage{hyperref}
\usepackage[nameinlink,capitalize,noabbrev]{cleveref}
\hypersetup{colorlinks=true,linkcolor=blue,citecolor=blue,urlcolor=blue}

\usepackage{tikz}
\usetikzlibrary{arrows.meta,calc}

\newcommand{\Rz}{\mathcal{R}_0}
\newcommand{\Rzh}{\mathcal{R}_0^{(h)}}
\newcommand{\Rze}{\mathcal{R}_0^{(e)}}
\newcommand{\Rtil}{\widetilde{\mathcal{R}}_0}
\newcommand{\Reff}{\mathcal{R}_{\mathrm{eff}}}
\newcommand{\dd}{\mathrm{d}}

\usepackage{authblk}

\title{Threshold Dynamics and Integrated Control of Cholera Transmission with Sanitation and Vaccination}
	
\author[1]{Abdallah Alsammani\thanks{Corresponding author: \texttt{aalsammani@desu.edu}}}
\author[2,3]{Gassan A.M.O. Farah}
\author[4]{Mohammed A.Y. Mohammed}
\author[5,6]{Mehmet Yavuz}

\affil[1]{Division of Physics, Engineering, Mathematics, and Computer Science, Delaware State University, Dover, DE, USA}
\affil[2]{Department of Mathematics and Applied Mathematics, University of the Western Cape, South Africa}
\affil[3]{National Center for Research, Institute of Space Research and Aerospace, Khartoum, Sudan}
\affil[4]{Department of Mathematics and Statistics, Georgia State University, Atlanta, GA, USA}
\affil[5]{Department of Mathematics and Computer Sciences, Necmettin Erbakan University, Konya 42090, T\"{u}rkiye}
\affil[6]{Department of Applied Mathematics and Informatics, Kyrgyz-Turkish Manas University, Bishkek 720038, Kyrgyzstan}

\date{}

\begin{document}
\maketitle

\begin{abstract}
\noindent
Cholera transmission is sustained by interacting human and environmental pathways, making integrated control essential in settings with inadequate water, sanitation, and hygiene infrastructure. We develop and analyze a Susceptible--Infected--Recovered--Water (SIWR) model with Monod-type environmental transmission and three interventions: human sanitation, environmental sanitation, and vaccination. The control reproduction number decomposes as $\Rz=\Rzh+\Rze$, separating direct and environmental transmission and yielding closed-form elasticity measures of pathway-specific contributions. We further quantify interactions between sanitation measures through an explicit Bliss independence synergy index. We establish local and global stability of the disease-free equilibrium and prove that the model undergoes a forward transcritical bifurcation at $\Rz=1$, excluding backward bifurcation. Global sensitivity analysis identifies direct transmission, the half-saturation constant, pathogen shedding, and human sanitation as major determinants of epidemic burden. Numerical simulations demonstrate that coordinated sanitation interventions can produce synergistic reductions in epidemic burden, providing a quantitative framework for evaluating integrated cholera control strategies.
\end{abstract}

\noindent\textbf{Keywords:} Cholera dynamics; Environmental transmission; Reproduction number; Bifurcation analysis; Integrated control

\section{Introduction}\label{sec1}


Cholera, an acute diarrhoeal disease caused by the bacterium \textit{Vibrio cholerae}, remains a major public health challenge in regions with inadequate water, sanitation, and hygiene (WASH) infrastructure. The World Health Organization estimates that $1.3$--$4$ million cases and $21{,}000$--$143{,}000$ deaths occur annually worldwide \cite{who2017cholera}. The burden is concentrated in sub-Saharan Africa, South Asia, and the Caribbean. At the same time, recent epidemics in Yemen, Haiti, and Zimbabwe illustrate the potential for rapid spread during conflict, population displacement, and disruption of essential public health services \cite{morris2010cholera,nelson2009transmission}.

Cholera spreads through two mechanistically distinct pathways: direct human transmission, including household and food-mediated contamination, and indirect transmission through contaminated aquatic reservoirs. These pathways operate on different timescales and respond differently to interventions because environmental reservoirs accumulate shed pathogens and sustain population-level exposure \cite{codeco2001endemic,eisenberg2013environmental}. Consequently, intervention effectiveness depends on the relative contribution of each pathway, which can vary substantially across epidemiological settings. Human transmission, for example, dominated the 2008--2009 Zimbabwe epidemic across all ten provinces \cite{mukandavire2011estimating}, whereas environmental transmission can predominate in estuarine and coastal settings \cite{codeco2001endemic}.

Mathematical models have been central to separating these transmission mechanisms and evaluating cholera control strategies. Early models introduced saturating environmental transmission and showed that aquatic reservoirs can sustain pathogen persistence even in the absence of direct transmission \cite{capasso1979mathematical,codeco2001endemic}. Subsequent SIWR models demonstrated that the relative contributions of direct and environmental transmission strongly influence intervention effectiveness \cite{tien2010multiple}, while outbreak studies estimated pathway-specific reproduction numbers from epidemiological data \cite{mukandavire2011estimating}. Other work has shown that asymptomatic infection complicates reproduction-number estimation \cite{king2008inapparent}, and vaccination has been investigated through simulation-based policy analyses \cite{azman2015impact}. Despite these advances, several analytical questions concerning transmission pathways and combined interventions remain unresolved.

First, although the additive decomposition $\Rz=\Rzh+\Rze$ separates human and environmental contributions to the control reproduction number, its implications for parameter sensitivity and pathway-specific intervention effects have received limited analytical treatment. Second, global stability results for cholera models with vaccination are often obtained under conditions stronger than $\Rz\le1$, raising the question of whether this gap reflects the Lyapunov construction or the underlying dynamics. Third, synergistic effects of combined interventions are commonly inferred from comparisons of epidemic trajectories rather than quantified through an explicit measure of interaction. Addressing these questions is important for connecting epidemic thresholds with the design and interpretation of integrated control strategies.

In this study, we analyze a parsimonious four-compartment SIWR model incorporating human sanitation, environmental sanitation, and vaccination. We derive closed-form elasticities of $\Rz$ and show that the sensitivity to each transmission parameter is determined by its corresponding pathway contribution. We establish a sharpened global stability condition for the disease-free equilibrium under $\Rtil\le1$ and show that this threshold cannot be improved within the class of linear Lyapunov functions. We further prove that the model undergoes a forward transcritical bifurcation at $\Rz=1$, excluding backward bifurcation within the proposed framework. Finally, using the Bliss independence criterion, we derive an explicit measure of sanitation synergy on the reproduction-number scale and numerically examine its epidemiological implications. These results connect pathway-specific transmission, epidemic thresholds, and intervention interactions within a unified analytical framework for cholera control.

\section{Model Formulation}\label{sec2}

\subsection{Model description and assumptions}\label{subsec:model_desc}

The human population at time $t$ is partitioned into susceptible individuals $S(t)$, infected (symptomatic, infectious) individuals $I(t)$ and recovered individuals $R(t)$, with $N(t)=S(t)+I(t)+R(t)$. The concentration of \textit{V.\ cholerae} in the shared aquatic reservoir is denoted $W(t)$, in cells per millilitre. The model rests on the following assumptions.

\begin{enumerate}
\item The population is demographically open: susceptibles are recruited at a constant rate $\Lambda$ and every human compartment experiences natural mortality at per-capita rate $\mu$.
\item Transmission occurs by direct human-to-human contact at rate $\beta_1$, and by ingestion of contaminated water at a rate that saturates with pathogen concentration.
\item Human sanitation, hygiene promotion, handwashing, safe food handling, reduces direct transmission by the factor $(1-\epsilon_h)$, $\epsilon_h\in[0,1]$.
\item Environmental sanitation, chlorination, water treatment and excreta management, reduces the environment-to-human route by the factor $(1-\epsilon_w)$, with $\epsilon_w\in[0,1]$.
\item Susceptibles are vaccinated at rate $\nu\ge0$ and move directly to $R$; vaccine protection is taken to be complete.
\item Infected individuals recover at rate $\gamma$, die from cholera at rate $\delta$, and shed pathogen into the reservoir at rate $\theta$.
\item Pathogen is removed from the reservoir at rate $\sigma$ by decay, dilution, and outflow.
\item Immunity acquired by recovery or vaccination does not wane over the timescales considered.
\item Only symptomatic infections are explicitly represented. Asymptomatic carriage is not included, and the implications of this modeling assumption are discussed in \cref{sec5}.
\end{enumerate}

Assumption 2 is made precise as follows.

\begin{definition}[Environment-to-human transmission function]\label{def:monod}
For $W\ge0$ the environment-to-human transmission rate is the Monod function
\begin{equation}\label{eq:monod}
\beta_2(W)=\beta_{\max}\frac{W}{k+W},
\end{equation}
where $\beta_{\max}>0$ is the maximal transmission rate and $k>0$ is the half-saturation concentration, in the form introduced by Monod for nutrient-limited bacterial growth \cite{Monod1949}. The function is smooth on $[0,\infty)$, strictly increasing, with $\beta_2(0)=0$, $\beta_2(k)=\beta_{\max}/2$, $\beta_2(W)\uparrow\beta_{\max}$, and
\begin{equation}\label{eq:monod_props}
\beta_2'(W)=\frac{\beta_{\max}k}{(k+W)^2}>0,\qquad \beta_2'(0)=\frac{\beta_{\max}}{k},\qquad \beta_2(W)\le \frac{\beta_{\max}}{k}\,W .
\end{equation}
\end{definition}

The last inequality in \eqref{eq:monod_props}, which states that the Monod function lies below its tangent at the origin, is used repeatedly in the stability analysis; it is exactly the concavity of $\beta_2$ that prevents the environmental route from generating the multiple endemic states associated with backward bifurcation.

Under these assumptions, the model reads
\begin{equation}\label{M1}
\begin{aligned}
\frac{\dd S}{\dd t} &= \Lambda - \beta_1(1-\epsilon_h)\frac{IS}{N} - \beta_2(W)(1-\epsilon_w)S - (\mu + \nu)S,\\
\frac{\dd I}{\dd t} &= \beta_1(1-\epsilon_h)\frac{IS}{N} + \beta_2(W)(1-\epsilon_w)S - (\gamma + \mu + \delta)I, \\
\frac{\dd R}{\dd t} &= \gamma I + \nu S - \mu R,\\
\frac{\dd W}{\dd t} &= \theta I - \sigma W,
\end{aligned}
\end{equation}
with $S(0)>0$, $I(0)\ge0$, $R(0)\ge0$, $W(0)\ge0$. Throughout we write
\begin{equation}\label{eq:Adef}
A:=\gamma+\mu+\delta
\end{equation}
for the total removal rate from the infected class, so that $1/A$ is the mean infectious period and $\delta/A$ is the case fatality proportion. The flow diagram is given in \cref{flowchart}.

\begin{figure}[!htbp]
\centering
\begin{tikzpicture}[>={Stealth[round]}, thick]
\tikzstyle{cmp} = [circle, minimum size=3.6em, text centered, font=\bfseries, thick, draw]

\node[cmp, fill=cyan!25]    (S) at (0,0)    {$S$};
\node[cmp, fill=red!30]     (I) at (5.0,0)  {$I$};
\node[cmp, fill=green!25]   (R) at (10.0,0) {$R$};
\node[cmp, fill=magenta!20] (W) at (5.0,-3.4) {$W$};

\draw[->] (-2.2,0) -- node[above] {$\Lambda$} (S);
\draw[->] (S) -- node[above] {$\beta_1(1-\epsilon_h)\dfrac{IS}{N}$} (I);
\draw[->] (I) -- node[above] {$\gamma I$} (R);
\draw[->] (S) to[bend left=38] node[above] {$\nu S$} (R);
\draw[->] (I) -- node[right] {$\theta I$} (W);
\draw[->] (W) to[bend right=22] node[below left=2pt] {$\beta_2(W)(1-\epsilon_w)S$} (S);
\draw[->] (W) -- node[above] {$\sigma W$} (8.4,-3.4);
\draw[->] (S) -- node[left] {$\mu S$} (0,-2.0);
\draw[->] (I) -- node[right] {$(\mu+\delta)I$} (6.6,-1.6);
\draw[->] (R) -- node[right] {$\mu R$} (10.0,-2.0);
\end{tikzpicture}
\caption{Flow diagram of model \eqref{M1}. Susceptible individuals acquire infection along two routes: direct contact, at rate $\beta_1(1-\epsilon_h)IS/N$, and ingestion of contaminated water, at rate $\beta_2(W)(1-\epsilon_w)S$ with $\beta_2$ given by \eqref{eq:monod}. Infected individuals shed pathogen into the reservoir at rate $\theta$, and pathogen is removed at rate $\sigma$. The three controls enter at distinct points: $\epsilon_h$ on the contact route, $\epsilon_w$ on the waterborne route, and $\nu$ on the susceptible pool, which is the structural reason why the first two controls interact while vaccination acts multiplicatively on both (Proposition \ref{prop:synergy}).}
\label{flowchart}
\end{figure}

\subsection{Parameter estimation and interpretation}\label{subsec:params}

Parameter ranges are collected in \cref{T1}. Several of the results below depend on the internal consistency of these values, so we state the reasoning explicitly.

\textbf{Demography ($\Lambda$, $\mu$).} These parameters are not independent: they determine the demographic equilibrium $N^\star=\Lambda/\mu$, which is the size of the community being modeled and is an observable of the setting rather than a free parameter. We therefore fix the community size and set $\Lambda=\mu N^\star$. Taking a life expectancy of $60$ years gives $\mu=1/(60\times365)=4.57\times10^{-5}$~day$^{-1}$, and a community of $N^\star=10^4$ gives $\Lambda=0.457$~day$^{-1}$. Since $\mu^{-1}\approx2.2\times10^4$~days greatly exceeds the $100$-day horizon of the simulations, demography influences the endemic state but is essentially inert during an outbreak; this is confirmed by the sensitivity analysis of \cref{subsec:sensitivity}, where $\mu$ has no detectable effect on the attack ratio.

\textbf{Disease-induced mortality $\delta$ and recovery rate $\gamma$.} Untreated severe cholera has a case fatality ratio of $25$--$50\%$, falling below $1\%$ with prompt rehydration \cite{harris2012cholera}. We use $\delta\in[0.01,0.05]$~day$^{-1}$, corresponding through $\delta/A$ to case fatality proportions of $4$--$15\%$ at $\gamma=0.25$~day$^{-1}$, which brackets the $5.3\%$ observed in Zimbabwe in 2008--2009. The symptomatic period of $2$--$10$ days gives $\gamma\in[0.1,0.5]$~day$^{-1}$ \cite{codeco2001endemic,azman2013incubation}.

\textbf{Transmission rates $\beta_1$ and $\beta_{\max}$.} Direct transmission covers fecal--oral routes not mediated by the shared reservoir, principally within households; outbreak-based estimates support $\beta_1\in[0.01,1.0]$~day$^{-1}$ \cite{tien2010multiple,mukandavire2011estimating}. The parameter $\beta_{\max}$ is the per-capita rate of infection from water use when the reservoir is saturated, and $\beta_{\max}\in[0.1,2.0]$~day$^{-1}$ spans settings from piped supply to exclusive reliance on surface water \cite{capasso1979mathematical}.

\textbf{Reservoir parameters $k$, $\theta$, $\sigma$.} The half-saturation concentration is tied to the infectious dose, which spans $10^3$--$10^8$ cells~mL$^{-1}$ depending on gastric acidity, serogroup and vehicle; we use $k\in[10^5,10^7]$ cells~mL$^{-1}$ following \cite{codeco2001endemic}. Symptomatic patients excrete $10^7$--$10^9$ cells~mL$^{-1}$ of stool at volumes of $10$--$20$~L~day$^{-1}$, and after dilution into a shared reservoir the effective contribution per patient is $\theta\in[1,100]$ cells~mL$^{-1}$~day$^{-1}$ \cite{nelson2009transmission}. Environmental survival ranges from hours in chlorinated water to weeks in warm nutrient-rich water, giving $\sigma\in[0.1,1.0]$~day$^{-1}$.

\textbf{Controls $\epsilon_h$, $\epsilon_w$, $\nu$.} The sanitation efficacies are dimensionless multipliers on the respective transmission terms and are varied over $[0,1]$. For the vaccination rate, we use $\nu\in[0,0.03]$~day$^{-1}$: the upper value vaccinates $95\%$ of the remaining susceptible pool over $100$ days and represents a well-resourced reactive campaign, whereas rates an order of magnitude larger, although sometimes assumed in modeling studies, exceed documented oral cholera vaccine delivery capacity \cite{azman2015impact,who2017cholera}.

\begin{table}[!htbp]
\centering
\caption{State variables and parameters of model \eqref{M1}. The recruitment rate is determined by $\Lambda=\mu N^\star$, where $N^\star$ is the community size, ensuring that the demographic equilibrium is consistent with the modeled population.}\label{T1}
\small
\begin{tabular}{@{}llll@{}}
\toprule
\textbf{Variable} & \textbf{Description} & \textbf{Unit} & \\
\midrule
$S(t)$ & Susceptible individuals & individuals & \\
$I(t)$ & Symptomatic infected individuals & individuals & \\
$R(t)$ & Recovered and vaccinated individuals & individuals & \\
$W(t)$ & Pathogen concentration in the reservoir & cells/mL & \\
\midrule
\textbf{Parameter} & \textbf{Description} & \textbf{Unit} & \textbf{Range} \\
\midrule
$N^\star$ & Community size, $\Lambda/\mu$ & individuals & $10^4$ \\
$\mu$ & Natural mortality rate & day$^{-1}$ & $[10^{-5},10^{-3}]$ \cite{codeco2001endemic} \\
$\delta$ & Disease-induced mortality rate & day$^{-1}$ & $[0.01,0.05]$ \cite{harris2012cholera} \\
$\gamma$ & Recovery rate & day$^{-1}$ & $[0.1,0.5]$ \cite{codeco2001endemic,azman2013incubation} \\
$\beta_1$ & Direct transmission rate & day$^{-1}$ & $[0.01,1.0]$ \cite{tien2010multiple,mukandavire2011estimating} \\
$\beta_{\max}$ & Maximal environmental transmission rate & day$^{-1}$ & $[0.1,2.0]$ \cite{capasso1979mathematical} \\
$k$ & Half-saturation concentration & cells/mL & $[10^5,10^7]$ \cite{codeco2001endemic} \\
$\theta$ & Shedding rate per infected individual & cells/(mL$\cdot$day) & $[1,100]$ \cite{nelson2009transmission} \\
$\sigma$ & Pathogen removal rate & day$^{-1}$ & $[0.1,1.0]$ \cite{codeco2001endemic} \\
$\epsilon_h$ & Human sanitation efficacy & -- & $[0,1]$ \\
$\epsilon_w$ & Environmental sanitation efficacy & -- & $[0,1]$ \\
$\nu$ & Vaccination rate & day$^{-1}$ & $[0,0.03]$ \cite{azman2015impact} \\
\bottomrule
\end{tabular}
\end{table}

\subsection{Well-posedness and the absorbing region}\label{sec:existence}

Define
\begin{equation}\label{eq:Omega}
\Omega=\Big\{(S,I,R,W)\in\mathbb{R}^4_{\ge0}\;:\;S+I+R\le \tfrac{\Lambda}{\mu},\quad W\le \tfrac{\theta\Lambda}{\mu\sigma}\Big\}.
\end{equation}

\begin{theorem}[Positive invariance]\label{thm:invariance}
The region $\Omega$ is positively invariant under the flow of \eqref{M1}, and every solution with nonnegative initial data enters $\Omega$ in finite time or approaches it asymptotically.
\end{theorem}

\begin{proof}
\emph{Nonnegativity.} On the face $S=0$ we have $\dd S/\dd t=\Lambda>0$; on $I=0$, $\dd I/\dd t=\beta_2(W)(1-\epsilon_w)S\ge0$; on $R=0$, $\dd R/\dd t=\gamma I+\nu S\ge0$; and on $W=0$, $\dd W/\dd t=\theta I\ge0$. Hence the vector field points into $\mathbb{R}^4_{\ge0}$ on its boundary and the nonnegative orthant is positively invariant. In particular $S(t)>0$ for $t\ge0$ whenever $S(0)>0$, so $N(t)>0$ and the quotient $IS/N$ is well defined.

\emph{Boundedness.} Adding the first three equations of \eqref{M1} gives $\dd N/\dd t=\Lambda-\mu N-\delta I\le \Lambda-\mu N$, so by the comparison principle $N(t)\le \Lambda/\mu+(N(0)-\Lambda/\mu)e^{-\mu t}$; thus $N(t)\le\Lambda/\mu$ for all $t\ge0$ if $N(0)\le\Lambda/\mu$, and $N(t)\to$ the interval $[0,\Lambda/\mu]$ otherwise. Since $I\le N$, the last equation gives $\dd W/\dd t\le \theta\Lambda/\mu-\sigma W$, whence $W(t)\le\theta\Lambda/(\mu\sigma)$ for $W(0)\le\theta\Lambda/(\mu\sigma)$ and $\limsup_{t\to\infty}W(t)\le\theta\Lambda/(\mu\sigma)$ in general.
\end{proof}

The following refinement is what makes the sharpened global stability result of \cref{subsec:global} possible. It states that, whatever the initial susceptible load, the flow drives $S$ below the disease-free level $\Lambda/(\mu+\nu)$, which is strictly smaller than $\Lambda/\mu$ when vaccination is active.

\begin{lemma}[Absorbing region]\label{lem:absorbing}
The set
\begin{equation}\label{eq:Omega1}
\Omega^{\flat}=\Omega\cap\Big\{S\le \tfrac{\Lambda}{\mu+\nu}\Big\}
\end{equation}
is positively invariant, and every solution of \eqref{M1} starting in $\Omega$ satisfies $\limsup_{t\to\infty}S(t)\le\Lambda/(\mu+\nu)$. Consequently $\Omega^\flat$ attracts all solutions in $\Omega$, and it suffices to study the asymptotic dynamics on $\Omega^\flat$.
\end{lemma}

\begin{proof}
From the first equation of \eqref{M1}, $\dd S/\dd t\le\Lambda-(\mu+\nu)S$, so $S(t)\le \Lambda/(\mu+\nu)+(S(0)-\Lambda/(\mu+\nu))e^{-(\mu+\nu)t}$, which yields the asymptotic bound. On the face $S=\Lambda/(\mu+\nu)$ within $\Omega$ we have $\dd S/\dd t=-\beta_1(1-\epsilon_h)IS/N-\beta_2(W)(1-\epsilon_w)S\le0$, so $\Omega^\flat$ is positively invariant.
\end{proof}

\begin{theorem}[Global existence and uniqueness]\label{thm:global}
For every initial condition in $\Omega$ with $S(0)>0$, system \eqref{M1} has a unique solution that exists for all $t\ge0$ and remains in $\Omega$.
\end{theorem}

\begin{proof}
The right-hand side of \eqref{M1} is continuously differentiable on the open set $\{N>0\}\supset\{S>0\}$, the only possible singularity being the division by $N$, which is excluded by \cref{thm:invariance}. Local existence and uniqueness follow from the Picard--Lindel\"of theorem \cite{Coddington1955}, and the uniform a priori bounds of \cref{thm:invariance} exclude finite-time blow-up, so the local solution continues to $[0,\infty)$.
\end{proof}

\section{Mathematical Analysis}\label{sec3}

\subsection{Equilibria}\label{subsec:equilibria}

Setting the right-hand side of \eqref{M1} to zero with $I=W=0$ gives the unique disease-free equilibrium (DFE)
\begin{equation}\label{eq:DFE}
E_0=(S^\star,0,R^\star,0)=\left(\frac{\Lambda}{\mu+\nu},\;0,\;\frac{\nu\Lambda}{\mu(\mu+\nu)},\;0\right),
\qquad S^\star+R^\star=N^\star=\frac{\Lambda}{\mu}.
\end{equation}
The disease-induced mortality rate $\delta$ is absent from $E_0$ because the term $\delta I$ vanishes there. The susceptible fraction at the DFE is
\begin{equation}\label{eq:susfrac}
\frac{S^\star}{N^\star}=\frac{\mu}{\mu+\nu},
\end{equation}
the proportion of the birth cohort that escapes vaccination before dying; it is attained on the demographic timescale $1/(\mu+\nu)$, a point of some importance for the interpretation of $\Rz$ under vaccination (Remark \ref{rem:timescale}).

An endemic equilibrium $E^*=(S^*,I^*,R^*,W^*)$ with $I^*>0$ satisfies
\begin{equation}\label{eq:EEsystem}
W^*=\frac{\theta}{\sigma}I^*,\qquad
N^*=\frac{\Lambda-\delta I^*}{\mu},\qquad
S^*=\frac{\Lambda-A I^*}{\mu+\nu},\qquad
R^*=N^*-S^*-I^*,
\end{equation}
where the second and third relations follow from $\dd N/\dd t=0$ and from adding the $S$ and $I$ equations. Substituting \eqref{eq:EEsystem} into the $I$ equation and cancelling the factor $I^*>0$ reduces the equilibrium problem to the scalar equation
\begin{equation}\label{eq:EEscalar}
\mathcal{G}(I^*):=\beta_1(1-\epsilon_h)\frac{S^*}{N^*}
+(1-\epsilon_w)\frac{\beta_{\max}\theta\,S^*}{\sigma k+\theta I^*}-A=0,
\qquad 0<I^*<\min\Big\{\frac{\Lambda}{A},\frac{\Lambda}{\delta}\Big\}.
\end{equation}
Both terms of $\mathcal{G}$ are strictly decreasing in $I^*$ on the admissible interval: the ratio $S^*/N^*$ decreases because $A>\delta$, and the second term decreases in both numerator and denominator. Hence $\mathcal{G}$ is strictly decreasing and \eqref{eq:EEscalar} has at most one root. Since $\mathcal{G}(0^+)=A(\Rz-1)$ with $\Rz$ given by \eqref{eq:R0} below, a positive root exists precisely when $\Rz>1$. This establishes the global existence and uniqueness of the endemic equilibrium, complementing the local bifurcation analysis in \cref{subsec:bifurcation}; the root is numerically located in \cref{subsec:bifurcation_num}.

\subsection{The control reproduction number and its additive structure}\label{subsec:R0}

We apply the next-generation method \cite{vandendriessche2002reproduction} to the infection subsystem $(I,W)$. Linearising at $E_0$ and splitting into new infections $F$ and transitions $V$ gives
\[
F=\begin{pmatrix} \beta_1(1-\epsilon_h)\dfrac{\mu}{\mu+\nu} & \dfrac{(1-\epsilon_w)\Lambda\beta_{\max}}{k(\mu+\nu)}\\[8pt] 0&0\end{pmatrix},
\qquad
V=\begin{pmatrix} A & 0\\ -\theta & \sigma\end{pmatrix},
\]
where the $(1,1)$ entry of $F$ uses the DFE susceptible fraction \eqref{eq:susfrac} and the $(1,2)$ entry uses $\beta_2'(0)=\beta_{\max}/k$ from \eqref{eq:monod_props}. Since
\[
V^{-1}=\frac{1}{\sigma A}\begin{pmatrix}\sigma&0\\ \theta& A\end{pmatrix},
\qquad
K_L=FV^{-1}=\begin{pmatrix}\dfrac{\beta_1(1-\epsilon_h)\mu}{(\mu+\nu)A}+\dfrac{(1-\epsilon_w)\Lambda\beta_{\max}\theta}{k\sigma(\mu+\nu)A} & \dfrac{(1-\epsilon_w)\Lambda\beta_{\max}}{k\sigma(\mu+\nu)}\\[8pt] 0&0\end{pmatrix}
\]
is upper triangular, its spectral radius is the $(1,1)$ entry:
\begin{equation}\label{eq:R0}
\Rz=\Rzh+\Rze,
\end{equation}
\begin{equation}\label{eq:R0h}
\Rzh=\frac{\beta_1(1-\epsilon_h)}{A}\cdot\frac{\mu}{\mu+\nu},
\qquad
\Rze=\frac{(1-\epsilon_w)\beta_{\max}\theta}{k\sigma A}\cdot\frac{\Lambda}{\mu+\nu}.
\end{equation}
Strictly, \eqref{eq:R0} is a \emph{control} reproduction number, since it depends on $\epsilon_h$, $\epsilon_w$ and $\nu$; it reduces to the basic reproduction number when all three vanish. We retain the symbol $\Rz$ throughout and state control values explicitly where relevant.

Each term has a direct mechanistic reading. $\Rzh$ counts the direct secondary infections produced by one case during its infectious period $1/A$, discounted by the susceptible fraction \eqref{eq:susfrac}. $\Rze$ traces the environmental loop: a case sheds $\theta$ cells~mL$^{-1}$ per day for $1/A$ days, each unit of concentration persists for $1/\sigma$ days, and while present it infects susceptibles at rate $(1-\epsilon_w)\beta_{\max}S^\star/k$. That $\Rze$ contains the susceptible \emph{number} $S^\star=\Lambda/(\mu+\nu)$ rather than a fraction reflects the fact that $W$ is a concentration in a shared reservoir: the environmental route is density dependent, so doubling the population drawing on a fixed water source doubles $\Rze$, whereas $\Rzh$ is unchanged.

The additive form \eqref{eq:R0} is usually quoted as an interpretive convenience. The next result shows that it determines the entire first-order sensitivity structure of $\Rz$. Write
\begin{equation}\label{eq:weights}
p_h=\frac{\Rzh}{\Rz},\qquad p_e=\frac{\Rze}{\Rz},\qquad p_h+p_e=1,
\end{equation}
for the pathway weights, and let $\Upsilon_q=\partial\log\Rz/\partial\log q$ denote the elasticity of $\Rz$ with respect to a parameter $q$.

\begin{proposition}[Elasticities of $\Rz$]\label{prop:elasticity}
For model \eqref{M1},
\begin{align}
&\Upsilon_{\beta_1}=p_h,\qquad
\Upsilon_{\beta_{\max}}=\Upsilon_{\theta}=\Upsilon_{\Lambda}=p_e,\qquad
\Upsilon_{k}=\Upsilon_{\sigma}=-p_e,\label{eq:elas1}\\[2pt]
&\Upsilon_{\gamma}=-\frac{\gamma}{A},\qquad
\Upsilon_{\delta}=-\frac{\delta}{A},\qquad
\Upsilon_{\nu}=-\frac{\nu}{\mu+\nu},\qquad
\Upsilon_{\mu}=p_h-\frac{\mu}{\mu+\nu}-\frac{\mu}{A},\label{eq:elas2}\\[2pt]
&\frac{\partial\log\Rz}{\partial\epsilon_h}=-\frac{p_h}{1-\epsilon_h},\qquad
\frac{\partial\log\Rz}{\partial\epsilon_w}=-\frac{p_e}{1-\epsilon_w}.\label{eq:elas3}
\end{align}
\end{proposition}

\begin{proof}
Each of $\beta_1$, $(1-\epsilon_h)$ appears only in $\Rzh$ and to the first power; each of $\beta_{\max},\theta,\Lambda,(1-\epsilon_w)$ appears only in $\Rze$ and to the first power; and $k,\sigma$ appear only in $\Rze$, with exponent $-1$. If a parameter $q$ occurs only in $\Rzh$ with exponent $m$ then $\partial\log\Rz/\partial\log q=(\partial\Rzh/\partial\log q)/\Rz=m\Rzh/\Rz=mp_h$, and analogously for $\Rze$; this gives \eqref{eq:elas1} and \eqref{eq:elas3}. The parameters $\gamma,\delta$ enter both components only through $A$, and $\nu$ only through $\mu+\nu$, so their elasticities are $-\gamma/A$, $-\delta/A$ and $-\nu/(\mu+\nu)$ by $p_h+p_e=1$. Finally $\mu$ occurs in $A$, in $\mu+\nu$, and additionally in the numerator of $\Rzh$, so
$\Upsilon_\mu=p_h\big(1-\tfrac{\mu}{\mu+\nu}-\tfrac{\mu}{A}\big)+p_e\big(-\tfrac{\mu}{\mu+\nu}-\tfrac{\mu}{A}\big)=p_h-\tfrac{\mu}{\mu+\nu}-\tfrac{\mu}{A}$.
\end{proof}
These identities establish a direct analytical connection between the pathway decomposition of the control reproduction number and the sensitivity rankings observed in the numerical analysis.

Three consequences are used below. First, the sign of the sensitivity of $\Rz$ to each parameter is fixed independently of parameter values; in particular $\Upsilon_k=-p_e<0$, so any correctly computed sensitivity index for the half-saturation constant must be negative. We use this as an external validation of the global sensitivity analysis in \cref{subsec:sensitivity}. Second, the magnitudes of the transmission sensitivities are not free but equal the pathway weights, so estimating $p_h$ and $p_e$ in a given setting immediately ranks all transmission parameters without further computation. Third, because $\Upsilon_{\beta_1}+\Upsilon_{\beta_{\max}}=1$, the two routes compete for sensitivity: a setting dominated by the environmental route is automatically insensitive to direct-contact parameters, and conversely.

\begin{remark}[Which reproduction number applies on which timescale]\label{rem:timescale}
Because $\nu$ enters \eqref{eq:R0h} only through the DFE susceptible fraction \eqref{eq:susfrac}, and that fraction is attained on the demographic timescale $1/(\mu+\nu)$, $\Rz$ measures the long-run consequence of a vaccination program sustained indefinitely. For $\nu\gg\mu$ it becomes very small, at $\nu=0.03$~day$^{-1}$ and $\mu=4.57\times10^{-5}$~day$^{-1}$ one has $\mu/(\mu+\nu)=1.5\times10^{-3}$, and it is then not a threshold for an outbreak developing over months in an initially susceptible population. The relevant outbreak-scale quantity replaces the DFE susceptible levels by the current ones,
\begin{equation}\label{eq:Reff}
\Reff(t)=\frac{\beta_1(1-\epsilon_h)}{A}\frac{S(t)}{N(t)}+\frac{(1-\epsilon_w)\beta_{\max}\theta}{k\sigma A}\,S(t),
\end{equation}
and satisfies $\Reff(0)=\Rz\big|_{\nu=0}$ in a fully susceptible population, whatever the value of $\nu$. A reactive campaign therefore does not lower the initial threshold; it accelerates the decline of $\Reff(t)$. Both quantities are correct, but they answer different questions, and conflating them leads to the erroneous conclusion that a campaign prevents an outbreak that in fact still occurs. This is quantified in \cref{subsec:vaccination}.
\end{remark}

\subsection{Local stability of the disease-free equilibrium}\label{subsec:local}

\begin{theorem}[Local asymptotic stability of the DFE]\label{thm:local_DFE}
The equilibrium $E_0$ is locally asymptotically stable if $\Rz<1$ and unstable if $\Rz>1$.
\end{theorem}

\begin{proof}
The Jacobian of \eqref{M1} at $E_0$ is
\[
J_{E_0}=\begin{pmatrix}
-(\mu+\nu) & -\beta_1(1-\epsilon_h)\dfrac{\mu}{\mu+\nu} & 0 & -\dfrac{(1-\epsilon_w)\Lambda\beta_{\max}}{k(\mu+\nu)}\\[8pt]
0 & \beta_1(1-\epsilon_h)\dfrac{\mu}{\mu+\nu}-A & 0 & \dfrac{(1-\epsilon_w)\Lambda\beta_{\max}}{k(\mu+\nu)}\\[8pt]
\nu & \gamma & -\mu & 0\\[4pt]
0 & \theta & 0 & -\sigma
\end{pmatrix}.
\]
Columns one and three show that $J_{E_0}$ is block triangular with respect to the splitting $(S,R)$, $(I,W)$. The $(S,R)$ block is triangular with eigenvalues $-(\mu+\nu)<0$ and $-\mu<0$, and the remaining spectrum is that of
\[
J_{\mathrm{inf}}=\begin{pmatrix}
\beta_1(1-\epsilon_h)\dfrac{\mu}{\mu+\nu}-A & \dfrac{(1-\epsilon_w)\Lambda\beta_{\max}}{k(\mu+\nu)}\\[8pt]
\theta & -\sigma
\end{pmatrix}.
\]
By \eqref{eq:R0h}, $\operatorname{tr}J_{\mathrm{inf}}=A(\Rzh-1)-\sigma$ and $\det J_{\mathrm{inf}}=\sigma A(1-\Rz)$. If $\Rz<1$ then $\Rzh\le\Rz<1$, so the trace is negative and the determinant positive, and the Routh--Hurwitz criterion for a $2\times2$ matrix gives two eigenvalues with negative real part. If $\Rz>1$ then $\det J_{\mathrm{inf}}<0$, so $J_{\mathrm{inf}}$ has one positive and one negative real eigenvalue and $E_0$ is a saddle.
\end{proof}

\subsection{Global stability of the disease-free equilibrium}\label{subsec:global}

Linear Lyapunov functions of the infection variables, standard for models of this type \cite{ShuaiDriessche2013,alsammani2023mathematical}, require an upper bound on the force of infection along all trajectories, and the bound obtained determines the threshold that can be proved. We first record the sharpest threshold obtainable in this way, then show that it cannot be improved within the class. Define
\begin{equation}\label{eq:Rtilde}
\Rtil=\frac{\beta_1(1-\epsilon_h)}{A}+\Rze=\Rzh\,\frac{\mu+\nu}{\mu}+\Rze .
\end{equation}
Thus $\Rtil$ is obtained from $\Rz$ by replacing the susceptible \emph{fraction} $\mu/(\mu+\nu)$ in the direct term by its supremum $1$ over $\Omega$, while retaining the sharp bound $S\le S^\star$ of Lemma \ref{lem:absorbing} in the environmental term. In particular $\Rz\le\Rtil$ always, with equality if and only if $\nu=0$.

\begin{theorem}[Global asymptotic stability of the DFE]\label{thm:global_DFE}
If $\Rtil\le1$ then $E_0$ is globally asymptotically stable in $\Omega$. When $\nu=0$ the condition reduces to $\Rz\le1$, which is sharp by \cref{thm:local_DFE}.
\end{theorem}

\begin{proof}
By Lemma \ref{lem:absorbing} it suffices to prove global attractivity on the absorbing set $\Omega^\flat$, on which $S\le S^\star=\Lambda/(\mu+\nu)$. Define on $\Omega^\flat$
\begin{equation}\label{eq:lyap}
V(I,W)=I+cW,\qquad c:=\frac{(1-\epsilon_w)\beta_{\max}S^\star}{k\sigma}=\frac{(1-\epsilon_w)\beta_{\max}\Lambda}{k\sigma(\mu+\nu)}>0 .
\end{equation}
Then $V$ is continuous and positive definite with respect to $\{I=W=0\}$. Differentiating along solutions of \eqref{M1},
\[
\dot V=\beta_1(1-\epsilon_h)\frac{IS}{N}+(1-\epsilon_w)\beta_2(W)S-AI+c(\theta I-\sigma W).
\]
Two bounds are applied. First, $S\le N$ gives $\beta_1(1-\epsilon_h)IS/N\le\beta_1(1-\epsilon_h)I$. Second, the tangent bound $\beta_2(W)\le\beta_{\max}W/k$ of \eqref{eq:monod_props} together with $S\le S^\star$ on $\Omega^\flat$ gives $(1-\epsilon_w)\beta_2(W)S\le(1-\epsilon_w)\beta_{\max}S^\star W/k=c\sigma W$. The two terms in $W$ cancel exactly, leaving
\[
\dot V\le\big[\beta_1(1-\epsilon_h)+c\theta-A\big]I
=A\Big[\frac{\beta_1(1-\epsilon_h)}{A}+\Rze-1\Big]I=A(\Rtil-1)I,
\]
where $c\theta=(1-\epsilon_w)\beta_{\max}\theta\Lambda/(k\sigma(\mu+\nu))=A\Rze$ by \eqref{eq:R0h}. Hence $\Rtil\le1$ implies $\dot V\le0$ on $\Omega^\flat$.

It remains to identify the largest invariant set in $\{\dot V=0\}$. If $\Rtil<1$, then $\dot V=0$ forces $I=0$. If $\Rtil=1$, then $\dot V=0$ forces equality in both bounds, that is $I(N-S)=0$ and $(1-\epsilon_w)S\big(\beta_{\max}W/k-\beta_2(W)\big)=0$; since $\beta_2$ is strictly concave the latter gives $W=0$, and then $\dot W=\theta I-\sigma W=0$ on an invariant set gives $I=0$ as well. In either case $I\equiv0$ on the invariant set, whence $\dot W=-\sigma W$ and $W\to0$, and the remaining equations reduce to $\dot S=\Lambda-(\mu+\nu)S$, $\dot R=\nu S-\mu R$, whose unique globally attracting equilibrium is $(S^\star,R^\star)$. The largest invariant set is therefore $\{E_0\}$, and LaSalle's invariance principle \cite{LaSalle1976} gives convergence of every solution in $\Omega^\flat$ to $E_0$. Together with local stability (\cref{thm:local_DFE}, applicable since $\Rz\le\Rtil\le1$) and Lemma \ref{lem:absorbing}, this yields global asymptotic stability in $\Omega$.
\end{proof}

The condition $\Rtil\le1$ is strictly weaker than the condition $\Rz(\mu+\nu)/\mu\le1$ that results from bounding $S$ by $N^\star=\Lambda/\mu$ in both terms, because the penalty factor $(\mu+\nu)/\mu$ is applied only to the direct-transmission component. The following proposition shows that the residual gap between $\Rz$ and $\Rtil$ is not an artifact of the particular function \eqref{eq:lyap}.

\begin{proposition}[Sharpness within the linear class]\label{prop:sharp}
Suppose $\Rtil>1$. Then for every $c\ge0$ the function $V_c=I+cW$ fails to be nonincreasing along the flow of \eqref{M1} in $\Omega^\flat$. Consequently, the threshold of \cref{thm:global_DFE} cannot be improved within the class of linear Lyapunov functions of $(I,W)$.
\end{proposition}

\begin{proof}
For $\varepsilon>0$ small consider $P_\varepsilon=(S^\star,\varepsilon,0,\theta\varepsilon/\sigma)$, which lies in $\Omega^\flat$ because $S^\star<\Lambda/\mu$. At $P_\varepsilon$ the reservoir is in quasi-equilibrium, $\dot W=\theta\varepsilon-\sigma(\theta\varepsilon/\sigma)=0$, so $\dot V_c=\dot I$ for every $c$. Since $N=S^\star+\varepsilon$ there,
\[
\frac{\dot I}{\varepsilon}\bigg|_{P_\varepsilon}
=\beta_1(1-\epsilon_h)\frac{S^\star}{S^\star+\varepsilon}
+(1-\epsilon_w)\frac{\beta_{\max}\theta S^\star}{\sigma k+\theta\varepsilon}-A
\;\xrightarrow[\varepsilon\to0^+]{}\;\beta_1(1-\epsilon_h)+A\Rze-A=A(\Rtil-1)>0,
\]
so $\dot V_c>0$ at $P_\varepsilon$ for all sufficiently small $\varepsilon>0$.
\end{proof}

\begin{remark}[Interpretation of the gap]\label{rem:gap}
Proposition \ref{prop:sharp} identifies the mechanism. Vaccination lowers $\Rz$ by depleting the susceptible \emph{fraction} at the DFE, but $\Omega$ contains states with $S/N$ close to one, in particular any population that has not yet been vaccinated, from which prevalence grows even when $\Rz<1$. The growth is transient and does not contradict global attractivity of $E_0$; it merely prevents a monotone linear functional of $(I,W)$ from certifying it. For instance, at $\nu=5\times10^{-3}$~day$^{-1}$ with the baseline values of \cref{tab:baseline} one has $\Rz=0.024$ and $\Rtil=1.86$; a solution seeded with ten cases in a fully susceptible community of $10^4$ produces an outbreak peaking at $1{,}765$ cases at day $22$ before converging to $E_0$, and along it $V$ increases by a factor of $177$ before decreasing to zero. Establishing global stability for all $\Rz<1$ therefore requires a nonlinear or composite Lyapunov construction, or a geometric approach; we do not claim such a result here.
\end{remark}

\subsection{Forward bifurcation at $\Rz=1$}\label{subsec:bifurcation}

We apply the center manifold theorem in the form of Castillo-Chavez and Song \cite{castilloChavezSong2004}, with bifurcation parameter $\phi=\beta_1$. Solving $\Rz=1$ for $\beta_1$ gives
\begin{equation}\label{eq:beta1star}
\beta_1^*=\frac{(1-\Rze)(\mu+\nu)A}{\mu(1-\epsilon_h)},
\end{equation}
which is positive precisely when $\Rze<1$, that is, when the environmental route alone cannot sustain transmission. Write $f=(f_S,f_I,f_R,f_W)$ for the right-hand side of \eqref{M1}.

\textbf{Eigenvectors.} At $(\beta_1^*,E_0)$ the matrix $J_{E_0}$ has a simple zero eigenvalue, all other eigenvalues having negative real part. From rows four, one and three of $J_{E_0}\mathbf{w}=\mathbf{0}$, a right null vector normalised by $w_I=1$ is
\begin{equation}\label{eq:rightvec}
w_I=1,\qquad w_W=\frac{\theta}{\sigma},\qquad
w_S=-\frac{A}{\mu+\nu}<0,\qquad
w_R=\frac{\gamma\mu-\nu(\mu+\delta)}{\mu(\mu+\nu)} .
\end{equation}
Columns one and three of $\mathbf{v}J_{E_0}=\mathbf{0}$ give $v_S=v_R=0$, and column four gives
\begin{equation}\label{eq:leftvec}
v_W=\frac{(1-\epsilon_w)\beta_{\max}S^\star}{k\sigma}\,v_I>0,
\end{equation}
with $v_I>0$ fixed by the normalisation $\mathbf{v}\cdot\mathbf{w}=1$. Note that $w_R$ changes sign: $w_R>0$ if and only if $\nu<\gamma\mu/(\mu+\delta)$, so at any realistic vaccination rate $w_R<0$. This matters because the second-order term carrying $w_R$ then contributes with a positive sign and cannot be discarded.

\textbf{Second derivatives.} Since $v_S=v_R=0$ and $f_W$ is linear, only $f_I$ contributes. Writing $B=\beta_1^*(1-\epsilon_h)$ and $C=(1-\epsilon_w)\beta_{\max}$, the nonvanishing second derivatives of $f_I$ at $E_0$ are
\begin{equation}\label{eq:secondderivs}
\begin{aligned}
\frac{\partial^2f_I}{\partial I^2}&=-\frac{2BS^\star}{(N^\star)^2}, &
\frac{\partial^2f_I}{\partial S\,\partial I}&=\frac{BR^\star}{(N^\star)^2}=\frac{B\nu}{N^\star(\mu+\nu)}, &
\frac{\partial^2f_I}{\partial R\,\partial I}&=-\frac{BS^\star}{(N^\star)^2},\\[2pt]
\frac{\partial^2f_I}{\partial W^2}&=-\frac{2CS^\star}{k^2}, &
\frac{\partial^2f_I}{\partial S\,\partial W}&=\frac{C}{k}. & &
\end{aligned}
\end{equation}
These follow from the identities
\[
\frac{\partial^2}{\partial I^2}\Big(\frac{IS}{N}\Big)=-\frac{2S(N-I)}{N^3},\qquad
\frac{\partial^2}{\partial S\,\partial I}\Big(\frac{IS}{N}\Big)=\frac{N-S}{N^2},\qquad
\frac{\partial^2}{\partial R\,\partial I}\Big(\frac{IS}{N}\Big)=-\frac{S}{N^2},
\]
evaluated at $I=0$, together with $\beta_2''(0)=-2\beta_{\max}/k^2$. The mixed derivative $\partial^2f_I/\partial I\partial W$ vanishes because $I$ and $W$ enter the two infection terms of $f_I$ separately.

\textbf{Bifurcation coefficients.} With $a=\sum_{k,i,j}v_kw_iw_j\,\partial^2f_k/\partial x_i\partial x_j$ and $b=\sum_{k,i}v_kw_i\,\partial^2f_k/\partial x_i\partial\phi$ evaluated at $(\beta_1^*,E_0)$,
\begin{equation}\label{eq:bcoef}
b=v_Iw_I\frac{\partial^2f_I}{\partial I\,\partial\beta_1}=v_I\,\frac{(1-\epsilon_h)\mu}{\mu+\nu}>0 ,
\end{equation}
and, substituting \eqref{eq:rightvec} and \eqref{eq:secondderivs},
\[
\frac{a}{v_I}=\frac{2B}{(N^\star)^2}\Big[-S^\star(1+w_R)+w_SR^\star\Big]
-\frac{2C\theta^2S^\star}{k^2\sigma^2}-\frac{2C\theta A}{k\sigma(\mu+\nu)} .
\]
The bracket simplifies. Using $R^\star/S^\star=\nu/\mu$ and $1+w_R=\big[\mu(\mu+\gamma)-\nu\delta\big]/\big(\mu(\mu+\nu)\big)$,
\[
-S^\star(1+w_R)+w_SR^\star
=S^\star\,\frac{-\mu(\mu+\gamma)+\nu\delta-A\nu}{\mu(\mu+\nu)}
=-S^\star\,\frac{\gamma+\mu}{\mu},
\]
since $\delta-A=-(\gamma+\mu)$ and the factor $(\mu+\nu)$ cancels. Hence, using $S^\star/(N^\star)^2=\mu^2/\big(\Lambda(\mu+\nu)\big)$,
\begin{equation}\label{eq:acoef}
a=-v_I\left[\frac{2\mu(\gamma+\mu)\beta_1^*(1-\epsilon_h)}{\Lambda(\mu+\nu)}
+\frac{2(1-\epsilon_w)\beta_{\max}\theta^2\Lambda}{k^2\sigma^2(\mu+\nu)}
+\frac{2(1-\epsilon_w)\beta_{\max}\theta A}{k\sigma(\mu+\nu)}\right]<0 .
\end{equation}
Every bracketed term is strictly positive for admissible parameters with $\epsilon_h,\epsilon_w\in[0,1)$, so $a<0$ holds unconditionally, in particular for all $\nu\ge0$; the potentially destabilizing contribution of $w_R<0$ is canceled exactly by susceptible depletion.

\begin{theorem}[Forward transcritical bifurcation]\label{thm:bifurcation}
System \eqref{M1} undergoes a forward (supercritical) transcritical bifurcation at $\Rz=1$. Precisely: for $\Rz\le1$ the DFE is locally asymptotically stable and no endemic equilibrium exists in $\Omega$; for $\Rz>1$ the DFE is unstable and a unique endemic equilibrium $E^*$ exists in the interior of $\Omega$, locally asymptotically stable for $\Rz$ in a right neighbourhood of unity; and the endemic prevalence increases continuously from zero as $\Rz$ crosses unity. In particular, the model exhibits neither backward bifurcation nor sub-threshold bistability.
\end{theorem}

\begin{proof}
Existence and uniqueness of $E^*$ for $\Rz>1$, and its nonexistence for $\Rz\le1$, follow from the strict monotonicity of $\mathcal{G}$ in \eqref{eq:EEscalar}. The stability statements follow from \cref{thm:local_DFE} and Theorem 4.1 of \cite{castilloChavezSong2004}: the zero eigenvalue of $J_{E_0}$ at $\beta_1=\beta_1^*$ is simple, all remaining eigenvalues have negative real part, and $a<0$, $b>0$ by \eqref{eq:acoef} and \eqref{eq:bcoef}. This is the sign pattern characterizing a forward bifurcation in which the bifurcating endemic branch is locally asymptotically stable for $\phi$ slightly above $\phi^*$.
\end{proof}

\begin{remark}\label{rem:noBB}
The structural reason for the absence of backward bifurcation is visible in \eqref{eq:secondderivs}: every nonlinearity of the infection term is stabilizing. The Monod function is concave, so $\partial^2f_I/\partial W^2<0$; standard incidence is concave in $I$ at fixed $S$; and the only potentially destabilizing contribution, through $w_R<0$, is removed by the identity above. Mechanisms known to generate backward bifurcation in cholera models, imperfect or waning vaccine protection, saturating treatment capacity, or a reservoir replenished independently of shedding, are absent here by construction. Epidemiologically, the elimination target is therefore unambiguous: no reduction beyond $\Rz<1$ is required, and no hysteresis has to be overcome.
\end{remark}

\subsection{Interaction between pathway-targeted controls}\label{subsec:synergy_theory}

A claim that interventions act synergistically is meaningful only relative to a stated null model of non-interaction. We adopt Bliss independence, the standard reference for interventions acting through distinct mechanisms: if control $i$ applied alone leaves a fraction $r_i$ of a burden measure, then independence predicts that the combination leaves $\prod_i r_i$. For a burden measure $Y$ with uncontrolled value $Y_0$, define the prevented fraction, the Bliss excess and the interaction index
\begin{equation}\label{eq:PF}
\mathcal{P}=1-\frac{Y}{Y_0},\qquad
\mathcal{E}=\mathcal{P}_{\mathrm{obs}}-\Big(1-\prod_i(1-\mathcal{P}_i)\Big),\qquad
\chi=\frac{\prod_i(1-\mathcal{P}_i)}{1-\mathcal{P}_{\mathrm{obs}}},
\end{equation}
so that $\mathcal{E}>0$ and $\chi>1$ indicate synergy. On the $\Rz$ scale, the interaction is available in closed form.

\begin{proposition}[Exact Bliss excess for $\Rz$]\label{prop:synergy}
Let $p_h,p_e$ be the pathway weights \eqref{eq:weights} in the absence of control and take $Y=\Rz$. Then human and environmental sanitation interact synergistically, with excess
\begin{equation}\label{eq:blissR0}
\mathcal{E}_{\Rz}(\epsilon_h,\epsilon_w,\nu)=p_hp_e\,\epsilon_h\epsilon_w\cdot\frac{\mu}{\mu+\nu},
\end{equation}
whereas vaccination is exactly Bliss-independent of both sanitation measures. The scaled excess $p_hp_e\epsilon_h\epsilon_w$ is maximized over pathway compositions at $p_h=p_e=\tfrac12$, and vanishes if and only if one route is absent or one sanitation efficacy is zero.
\end{proposition}

\begin{proof}
By \eqref{eq:R0}--\eqref{eq:R0h} the fraction of the uncontrolled $\Rz$ remaining under joint control is
\[
r(\epsilon_h,\epsilon_w,\nu)=\big[p_h(1-\epsilon_h)+p_e(1-\epsilon_w)\big]\frac{\mu}{\mu+\nu}
=(1-p_h\epsilon_h-p_e\epsilon_w)\frac{\mu}{\mu+\nu},
\]
while the single-control values are $r_h=1-p_h\epsilon_h$, $r_w=1-p_e\epsilon_w$ and $r_\nu=\mu/(\mu+\nu)$. Hence
\[
\mathcal{E}_{\Rz}=r_hr_wr_\nu-r=\big[(1-p_h\epsilon_h)(1-p_e\epsilon_w)-(1-p_h\epsilon_h-p_e\epsilon_w)\big]\frac{\mu}{\mu+\nu}
=p_hp_e\epsilon_h\epsilon_w\frac{\mu}{\mu+\nu}.
\]
The factorization of $r$ into a sanitation factor times $r_\nu$ shows that vaccination contributes no interaction term. Finally $p_hp_e=p_h(1-p_h)\le\tfrac14$, with equality at $p_h=\tfrac12$.
\end{proof}

Proposition \ref{prop:synergy} makes the mechanism of synergy explicit and is, in our view, the most useful practical consequence of the additive decomposition. Because $\Rz$ is a \emph{sum} over routes, a control that blocks one route leaves the other intact, and its marginal return saturates: no amount of environmental sanitation can bring $\Rz$ below $\Rzh$. Independence, being multiplicative, implicitly credits each control with a proportional reduction of the whole; partial control of both routes therefore beats that expectation by exactly $p_hp_e\epsilon_h\epsilon_w$. Synergy is largest where the two routes contribute equally to risk and negligible where one dominates, a testable statement, since $p_h$ and $p_e$ are estimable from outbreak data \cite{mukandavire2011estimating}.

Two qualifications follow. First, \eqref{eq:blissR0} describes interaction on the threshold scale; interventions may exhibit stronger interactions on epidemic outcome scales such as attack ratio or peak prevalence because these outcomes depend nonlinearly on the reproduction number. Consequently, comparable reductions in transmission can produce disproportionately large reductions in epidemic burden as the system approaches the epidemic threshold. This additional dynamic contribution to synergy is quantified in \cref{subsec:combined}. Second, $\mathcal{E}_{\Rz}\le\tfrac14\epsilon_h\epsilon_w\le\tfrac14$, so the interaction on the threshold scale is real but bounded; larger observed synergies on epidemic outcome scales may therefore reflect the nonlinear relationship between transmission and epidemic burden in addition to the structure of $\Rz$.

\section{Numerical Results}\label{sec4}

\subsection{Implementation and baseline scenario}\label{subsec:numerics}

System \eqref{M1} was integrated over a $100$-day horizon with an adaptive stiff/non-stiff solver at relative tolerance $10^{-10}$ and absolute tolerance $10^{-12}$; halving the tolerances changed all reported quantities by less than $0.01\%$. Endemic equilibria were obtained by bracketing the root of the scalar equation \eqref{eq:EEscalar} and their stability was determined from the eigenvalues of the full Jacobian, rather than by long-time integration, which avoids the loss of accuracy that occurs when prevalence is small. Baseline values are listed in \cref{tab:baseline}.

\begin{table}[!htbp]
\centering
\caption{Baseline values used in the numerical experiments, representing an unserved urban community of $N^\star=10^4$ during an uncontrolled outbreak. The recruitment rate is $\Lambda=\mu N^\star$, so that the demographic equilibrium coincides with the modeled population and the initial condition lies in $\Omega$.}\label{tab:baseline}
\small
\begin{tabular}{@{}lll@{}}
\toprule
\textbf{Parameter} & \textbf{Value} & \textbf{Basis} \\
\midrule
$N^\star$ & $10^4$ individuals & Community served by one water source \\
$\mu$ & $4.57\times10^{-5}$ day$^{-1}$ & Life expectancy $60$ years \\
$\Lambda$ & $0.457$ day$^{-1}$ & $\mu N^\star$ \\
$\delta$ & $0.02$ day$^{-1}$ & Case fatality $\delta/A=7.4\%$ \\
$\gamma$ & $0.25$ day$^{-1}$ & $4$-day mean symptomatic period \\
$\beta_1$ & $0.50$ day$^{-1}$ & Moderate direct contact \\
$\beta_{\max}$ & $0.75$ day$^{-1}$ & Ceiling on waterborne transmission \\
$k$ & $10^6$ cells/mL & Half-saturation concentration \\
$\theta$ & $10$ cells/(mL$\cdot$day) & Effective shedding per case \\
$\sigma$ & $0.33$ day$^{-1}$ & $3$-day mean pathogen persistence \\
$\epsilon_h,\epsilon_w,\nu$ & $0$ & Uncontrolled baseline \\
\midrule
$S(0),I(0),R(0)$ & $9{,}990$, $10$, $0$ & Fully susceptible community, ten seeded cases \\
$W(0)$ & $100$ cells/mL & Low-level background contamination, $W(0)/k=10^{-4}$ \\
\bottomrule
\end{tabular}
\end{table}

At these values $\Rzh=1.85$ and $\Rze=0.84$, so $\Rz=2.69$ with pathway weights $p_h=0.69$ and $p_e=0.31$. Both the magnitude and the composition are consistent with the province-level estimates of Mukandavire et~al.\ \cite{mukandavire2011estimating} for the 2008--2009 Zimbabwean epidemic, where $\Rz$ ranged from just above $1$ to $2.72$ and human-to-human transmission was the larger contributor in a land-locked setting without an estuarine reservoir. The baseline therefore represents a community in which direct transmission dominates. Still, the environmental route is far from negligible, which is the regime in which the interaction of Proposition \ref{prop:synergy} is appreciable.

\textbf{Sensitivity to initial conditions.} Because the seeded prevalence and initial contamination are not observable quantities, we verified that they do not drive the conclusions. Varying $I(0)$ over $\{1,10,100\}$ and $W(0)$ over $\{0,10^2,10^4,10^5\}$ cells/mL changes the cumulative number of infections only between $9{,}098$ and $9{,}368$, a spread of under $3\%$, and the peak prevalence only between $2{,}171$ and $2{,}818$. What these quantities control is timing: the peak occurs on day $28.1$ for a single-seeded case with a clean reservoir, and on day $6.6$ for $100$ cases with $W(0)=10^5$ cells/mL. This separation is expected. In a fully susceptible population the epidemic is governed by the dominant eigenvalue of the linearised infection subsystem, here $0.341$~day$^{-1}$ (doubling time two days), so the initial data enter only through the offset of an exponential growth phase. At the same time, the final size is determined by $\Reff$ and the susceptible pool. Reported burdens are therefore robust, whereas statements about peak timing should be read as conditional on the seeding assumption.

\subsection{Baseline epidemic dynamics}\label{subsec:baseline}

\begin{figure}[!htbp]
\centering
\includegraphics[width=\textwidth]{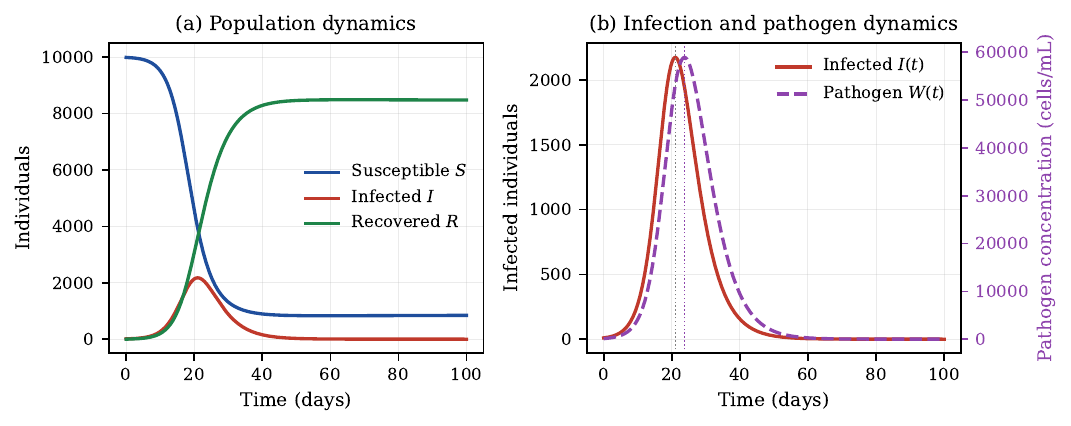}
\caption{Uncontrolled outbreak ($\epsilon_h=\epsilon_w=\nu=0$, $\Rz=2.69$), baseline values from \cref{tab:baseline}. \textbf{(a)} Susceptible, infected, and recovered individuals. Prevalence peaks at $2{,}175$ on day $21$ and the susceptible pool is depleted to $829$; the terminal population of $9{,}322$ reflects $678$ cholera deaths. \textbf{(b)} Prevalence (left axis) and reservoir concentration $W(t)$ (right axis); dotted vertical lines mark the two maxima. The reservoir peaks $2.7$ days after prevalence because \eqref{M1} makes $W$ a low-pass filter of $I$ with time constant $1/\sigma=3$ days, so contamination remains near its maximum while incidence is already declining.}
\label{fig:baseline}
\end{figure}

\Cref{fig:baseline} shows the uncontrolled outbreak. Prevalence peaks at $2{,}175$ cases on day $21$, that is, $21.8\%$ of the community is simultaneously infected, and the susceptible pool falls to $829$ individuals before demographic replenishment begins to matter. Over $100$ days $9{,}179$ infections occur, an attack ratio of $92\%$, and the population declines from $10{,}000$ to $9{,}322$, the $678$ deaths corresponding to the case fatality proportion $\delta/A=7.4\%$.

The pathogen concentration peaks at $5.9\times10^4$ cells~mL$^{-1}$ on day $23.8$, lagging prevalence by $2.7$ days. The lag is not a fitted feature but a direct consequence of the reservoir equation: $W$ is a low-pass filter of prevalence with time constant $1/\sigma=3$ days, so its peak occurs where $\theta I=\sigma W$, shortly after prevalence turns over. Two implications follow. First, water treatment retains value when deployed after the clinical peak, because the reservoir is still near its maximum while incidence is falling. Second, the peak concentration $5.9\times10^4$ cells~mL$^{-1}$ remains well below the half-saturation constant $k=10^6$, so $\beta_2$ operates in its near-linear regime throughout; the saturation nonlinearity therefore shapes the model's qualitative properties, it is what makes $\partial^2f_I/\partial W^2<0$ in \eqref{eq:secondderivs} and thereby excludes backward bifurcation, without materially distorting the baseline trajectory. Saturation would bind only in settings with a substantially more contaminated reservoir, larger $\theta$, or smaller $k$.

\subsection{Single interventions}\label{subsec:single}

\subsubsection{Human sanitation}\label{subsec:human_san}

\begin{figure}[!htbp]
\centering
\includegraphics[width=\textwidth]{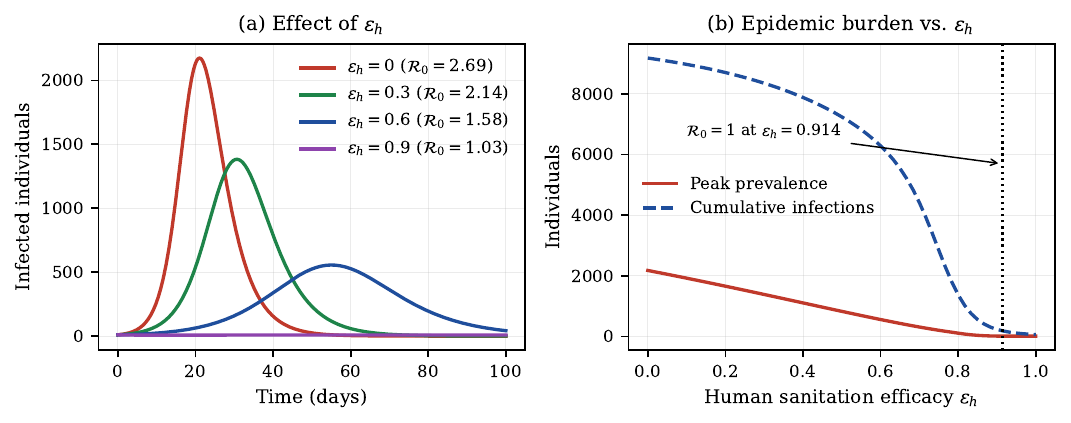}
\caption{Human sanitation. \textbf{(a)} Prevalence for $\epsilon_h=0,0.3,0.6,0.9$, all other parameters at baseline (\cref{tab:baseline}); the legend reports the corresponding $\Rz$. Increasing $\epsilon_h$ reduces the epidemic growth rate from $0.341$ to $0.004$~day$^{-1}$, thereby suppressing and delaying the peak. At $\epsilon_h=0.9$ no outbreak is visible over $100$ days despite $\Rz=1.03>1$, because the doubling time is then $168$ days. \textbf{(b)} Peak prevalence and cumulative infections against $\epsilon_h$. The dotted line marks $\epsilon_h=0.914$, above which $\Rz<1$; the threshold is high because human sanitation can reduce $\Rz$ only to its floor $\Rze=0.84$.}
\label{fig:human_sanitation}
\end{figure}

Raising $\epsilon_h$ from $0$ to $0.6$ lowers $\Rz$ from $2.69$ to $1.58$, reduces peak prevalence from $2{,}175$ to $556$ cases ($-74\%$) and cumulative infections from $9{,}179$ to $6{,}283$ ($-32\%$), while delaying the peak from day $21$ to day $55$ (\cref{fig:human_sanitation}). The disparity between the peak and cumulative reductions is important for interpretation: substantial epidemic flattening can occur before comparably large reductions in cumulative burden are achieved, reflecting the different dependence of peak prevalence and cumulative infections on transmission intensity over the finite simulation horizon.

The delay itself is governed by the dominant eigenvalue of the linearised infection subsystem, which falls from $0.341$~day$^{-1}$ at $\epsilon_h=0$ to $0.217$, $0.103$ and $0.0041$~day$^{-1}$ at $\epsilon_h=0.3,0.6,0.9$, corresponding to doubling times of $2$, $3$, $7$ and $168$ days. This explains an apparent anomaly in \cref{fig:human_sanitation}: at $\epsilon_h=0.9$ the epidemic does not take off at all within the horizon even though $\Rz=1.03>1$. The equilibrium theory is not violated; the endemic state is still approached, but with a doubling time of $168$ days, the approach is invisible over $100$ days. Threshold statements and horizon statements are distinct, and only the growth rate reconciles them.

By \eqref{eq:R0} human sanitation alone reduces $\Rz$ only to its floor $\Rze=0.84$, so elimination requires
\[
\epsilon_h>1-\frac{1-\Rze}{\Rzh}=0.914 ,
\]
an efficacy far beyond what hygiene promotion achieves in practice. Human sanitation is thus the more powerful of the two sanitation levers at this baseline, as Proposition \ref{prop:elasticity} predicts, since $|\partial\log\Rz/\partial\epsilon_h|=p_h=0.69$ exceeds $|\partial\log\Rz/\partial\epsilon_w|=p_e=0.31$. Yet, it cannot achieve elimination on its own.

\subsubsection{Environmental sanitation}\label{subsec:env_san}

\begin{figure}[!htbp]
\centering
\includegraphics[width=\textwidth]{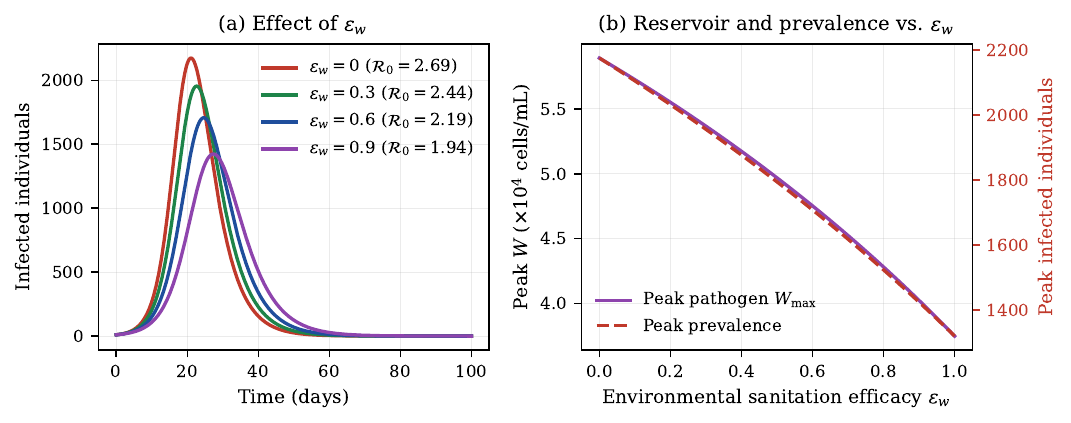}
\caption{Environmental sanitation. \textbf{(a)} Prevalence for $\epsilon_w=0,0.3,0.6,0.9$ with all other parameters at baseline (\cref{tab:baseline}); the legend reports $\Rz$. Because $\Rz\ge\Rzh=1.85$ for every $\epsilon_w$, this control cannot achieve elimination at the baseline pathway composition. \textbf{(b)} Peak reservoir concentration (left axis) and peak prevalence (right axis) against $\epsilon_w$. Contamination declines nonlinearly because $\epsilon_w$ acts on $W$ only indirectly, through reduced shedding by a smaller infected population.}
\label{fig:environmental_sanitation}
\end{figure}

Environmental sanitation at matched efficacy is markedly less effective at this baseline: $\epsilon_w=0.9$ lowers $\Rz$ only from $2.69$ to $1.94$, peak prevalence to $1{,}424$ cases ($-35\%$) and cumulative infections to $7{,}940$ ($-14\%$) (\cref{fig:environmental_sanitation}). Since $\Rz\ge\Rzh=1.85>1$ irrespective of $\epsilon_w$, environmental sanitation alone cannot eliminate transmission here at any efficacy. This is the mirror image of the previous subsection and is entirely determined by the pathway weights: whichever route carries the larger share of $\Rz$ defines the achievable floor for the control targeting the other route. Thus, the relative effectiveness of human and environmental sanitation is not a universal property of the model but depends on the pathway composition of the epidemiological setting. In estuarine or strongly reservoir-driven settings, where $p_e>p_h$, the ranking may reverse.

The reservoir responds more strongly than prevalence: peak concentration falls from $5.9\times10^4$ to $4.0\times10^4$ cells~mL$^{-1}$ as $\epsilon_w$ rises to $0.9$, and the decline steepens at high efficacy. The mechanism is a feedback rather than a saturation effect: $\epsilon_w$ does not enter the $W$ equation at all, so contamination falls only because fewer infected individuals shed. This indirect route strengthens as prevalence itself begins to respond. The practical corollary is that measurements of reservoir contamination understate the epidemiological benefit of water treatment at low efficacy and overstate its marginal benefit at high efficacy.

\subsubsection{Vaccination}\label{subsec:vaccination}

\begin{figure}[!htbp]
\centering
\includegraphics[width=\textwidth]{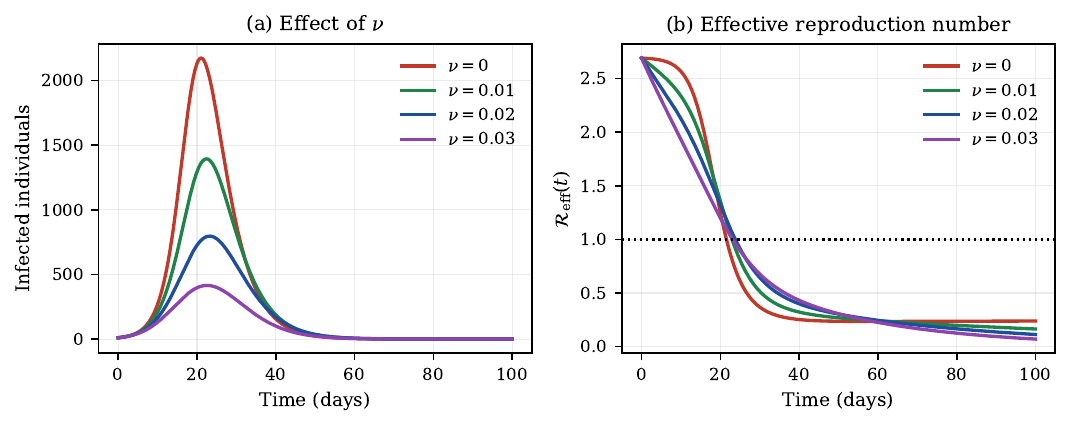}
\caption{Vaccination. \textbf{(a)} Prevalence for $\nu=0,0.01,0.02,0.03$~day$^{-1}$ with all other parameters at baseline (\cref{tab:baseline}). An outbreak occurs at every rate, although cumulative infections fall by up to $72\%$. \textbf{(b)} Effective reproduction number \eqref{eq:Reff}. All four curves start at $\Reff(0)=2.69$ because the population is initially fully susceptible; $\Reff$ crosses unity (dotted line) between days $21.5$ and $23.8$ in all cases, but at $\nu=0.03$ some $73\%$ of the immunity responsible for that decline is acquired by vaccination rather than by infection. The control reproduction numbers computed at the disease-free equilibrium, $\Rz\le0.012$ for $\nu\ge0.01$, refer to the vaccinated demographic steady state and are not outbreak thresholds (Remark \ref{rem:timescale}).}
\label{fig:vaccination}
\end{figure}

Vaccination at $\nu=0.01$, $0.02$ and $0.03$~day$^{-1}$ reduces peak prevalence by $36\%$, $63\%$ and $81\%$, and cumulative infections by $25\%$, $51\%$ and $72\%$ (\cref{fig:vaccination}). Critically, the outbreak still occurs at every tested rate, even though the corresponding control reproduction numbers computed from \eqref{eq:R0} are $0.012$, $0.006$, and $0.004$. This is exactly the situation anticipated in Remark \ref{rem:timescale}. Those values describe the vaccinated demographic steady state, which for $\mu^{-1}=60$ years lies decades away; they say nothing about a community that is fully susceptible when the campaign starts. \Cref{fig:vaccination}(b) makes the point directly: $\Reff(0)=2.69$ for all four rates, and $\Reff$ crosses unity at a comparable time in every case, between days $21.5$ and $23.8$. What changes with vaccination is not when the susceptible pool is depleted but \emph{how}. The reproduction number falls because susceptibles are removed, and in the uncontrolled outbreak the only removal mechanism is infection; vaccination substitutes for it. At $\nu=0.01$, $0.02$ and $0.03$~day$^{-1}$ the proportion of acquired immunity obtained by vaccination rather than by infection is $27\%$, $53\%$ and $73\%$ respectively, and cumulative infections fall in the same proportion.

This distinction has a concrete programmatic consequence. Because a campaign competes with the epidemic to deplete the same susceptible pool, its benefit depends on the interval between deployment and the epidemic peak: doses delivered after infection has already conferred immunity avert nothing. At the baseline growth rate of $0.341$~day$^{-1}$, that interval is short, which is the quantitative form of the familiar argument for pre-positioning oral cholera vaccine stocks rather than relying on reactive deployment \cite{azman2015impact,who2017cholera}. It also explains why vaccination, despite having the largest single-control effect on cumulative infections in \cref{tab:synergy}, is not the parameter most sensitive to the $100$-day attack ratio (\cref{subsec:sensitivity}).

\subsection{Combined interventions and quantification of synergy}\label{subsec:combined}

\begin{figure}[!htbp]
\centering
\includegraphics[width=0.7\textwidth]{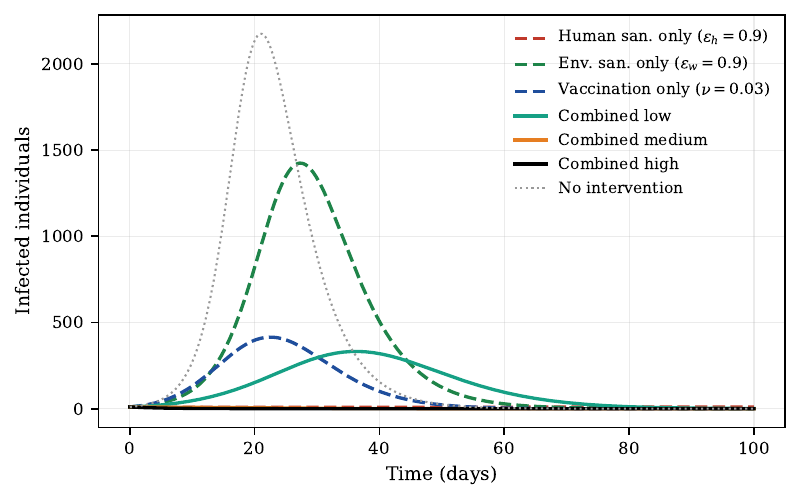}
\caption{Single and combined control strategies, baseline parameters otherwise as in \cref{tab:baseline}. Dashed curves show each control alone at its maximum tested intensity ($\epsilon_h=0.9$; $\epsilon_w=0.9$; $\nu=0.03$); solid curves show combinations at low ($0.3,0.3,0.01$), medium ($0.6,0.6,0.02$) and high ($0.9,0.9,0.03$) intensity; the dotted curve is the uncontrolled outbreak. Medium-intensity combined control extinguishes the outbreak with an efficacy of $0.6$, whereas environmental sanitation at $0.9$ leaves a peak of $1{,}424$ cases. Quantitative interaction metrics are given in \cref{tab:synergy}.}
\label{fig:combined_interventions}
\end{figure}

We compare each control applied alone with three combinations of increasing intensity: low ($\epsilon_h=\epsilon_w=0.3$, $\nu=0.01$), medium ($0.6$, $0.6$, $0.02$) and high ($0.9$, $0.9$, $0.03$). Trajectories are shown in \cref{fig:combined_interventions} and the interaction metrics \eqref{eq:PF}, computed on cumulative infections, are collected in \cref{tab:synergy}.

\begin{table}[!htbp]
\centering
\caption{Quantified interaction between the three controls, using cumulative infections over $100$ days as the burden measure and the Bliss independence null model \eqref{eq:PF}. $\mathcal{P}_h,\mathcal{P}_e,\mathcal{P}_\nu$ are prevented fractions for each control applied alone at the stated intensity; $\mathcal{P}_{\mathrm{obs}}$ is the observed prevented fraction of the combination and $\mathcal{P}_{\mathrm{Bliss}}=1-\prod_i(1-\mathcal{P}_i)$ the independence prediction. Positive excess $\mathcal{E}$ and interaction index $\chi>1$ indicate synergy.}\label{tab:synergy}
\small
\begin{tabular}{@{}lccccccc@{}}
\toprule
Intensity & $(\epsilon_h,\epsilon_w,\nu)$ & $\mathcal{P}_h$ & $\mathcal{P}_e$ & $\mathcal{P}_\nu$ & $\mathcal{P}_{\mathrm{Bliss}}$ & $\mathcal{P}_{\mathrm{obs}}$ & $\mathcal{E}$ \ ($\chi$) \\
\midrule
Low    & $(0.3,0.3,0.01)$ & $0.090$ & $0.030$ & $0.254$ & $0.342$ & $0.657$ & $0.315$ \ ($1.9$) \\
Medium & $(0.6,0.6,0.02)$ & $0.316$ & $0.073$ & $0.506$ & $0.686$ & $0.995$ & $0.309$ \ ($61$) \\
High   & $(0.9,0.9,0.03)$ & $0.975$ & $0.135$ & $0.717$ & $0.994$ & $1.000$ & $0.006$ \ ($17$) \\
\bottomrule
\end{tabular}
\end{table}

Three findings emerge. First, synergy is real and substantial at low and moderate intensity. At the low setting, the three controls used singly prevent $9.0\%$, $3.0\%$ and $25.4\%$ of infections, so independence predicts $34.2\%$ for the combination; the observed value is $65.7\%$, an excess of $31.5$ percentage points and an interaction index of $\chi=1.9$. At medium intensity, the combination prevents $99.5\%$ of infections and extinguishes the outbreak, whereas independence predicts $68.6\%$.

Second, the combination outperforms every single control applied at maximum intensity. Medium-intensity combined control prevents $99.5\%$ of infections using efficacies of $0.6$, while environmental sanitation at $\epsilon_w=0.9$ prevents only $13.5\%$ and vaccination at $\nu=0.03$ only $71.7\%$. Human sanitation at $\epsilon_h=0.9$ prevents $97.5\%$, so it is not the case that the combination dominates every alternative by a wide margin; but achieving $\epsilon_h=0.9$ means eliminating nine of every ten direct transmission events, which is not attainable by hygiene promotion in practice, whereas three simultaneous efficacies of $0.6$ are.

\begin{figure}[!htbp]
\centering
\includegraphics[width=\textwidth]{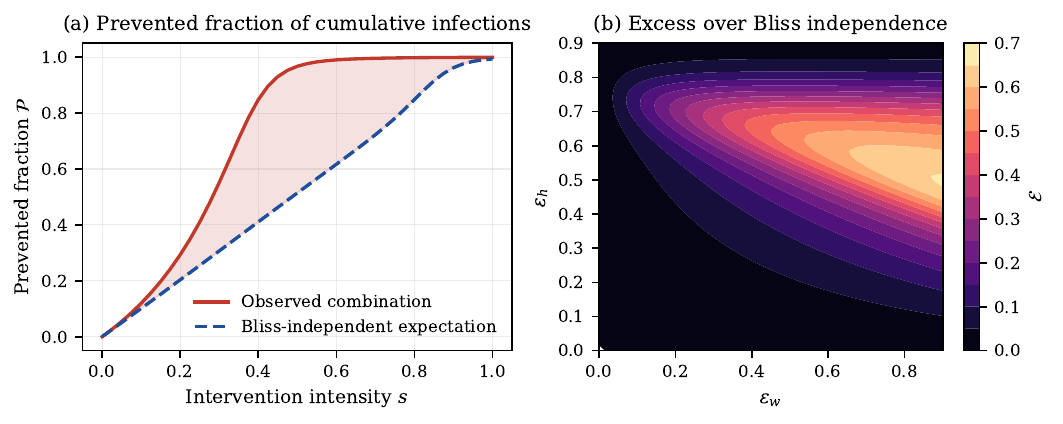}
\caption{Quantified synergy on the burden scale, using cumulative infections. \textbf{(a)} Prevented fraction of the joint strategy $(\epsilon_h,\epsilon_w,\nu)=(0.9s,0.9s,0.03s)$ against intensity $s$, compared with the Bliss independence prediction from the three single-control simulations; the shaded area is the excess $\mathcal{E}$. The gap opens from zero, peaks at intermediate intensity, and closes as both curves saturate near $\mathcal{P}=1$. \textbf{(b)} Excess $\mathcal{E}$ for the two sanitation controls over $(\epsilon_h,\epsilon_w)$ at $\nu=0$, reaching $0.65$. The ridge corresponds to combinations that cross the elimination threshold while neither control does so alone.}
\label{fig:synergy}
\end{figure}

Third, the interaction is not driven solely by the structure of $\Rz$. Proposition \ref{prop:synergy} gives, at the low setting, $p_hp_e\epsilon_h\epsilon_w=0.019$ on the $\Rz$ scale, substantially smaller than the observed excess of $0.315$ on the burden scale. This difference reflects the nonlinear relationship between transmission intensity and epidemic burden, through which combined reductions in transmission can produce disproportionately large reductions in cumulative infections. \Cref{fig:synergy}(a) traces this behavior over a continuous intensity path and shows the gap between the observed and independent prevented fractions opening from zero, reaching a maximum near intermediate intensity, and closing again as both curves approach saturation.

\Cref{fig:synergy}(b) shows the excess for the two sanitation controls over the $(\epsilon_h,\epsilon_w)$ plane at $\nu=0$, reaching $0.65$ near $\epsilon_h\approx0.5$ and $\epsilon_w\approx0.9$. The ridge lies where the combination crosses the elimination threshold. At the same time, neither control does so alone, which is precisely the region Proposition \ref{prop:synergy} identifies as the operational value of integrated control.

\subsection{Bifurcation structure}\label{subsec:bifurcation_num}

\begin{figure}[!ht]
\centering
\includegraphics[width=\textwidth]{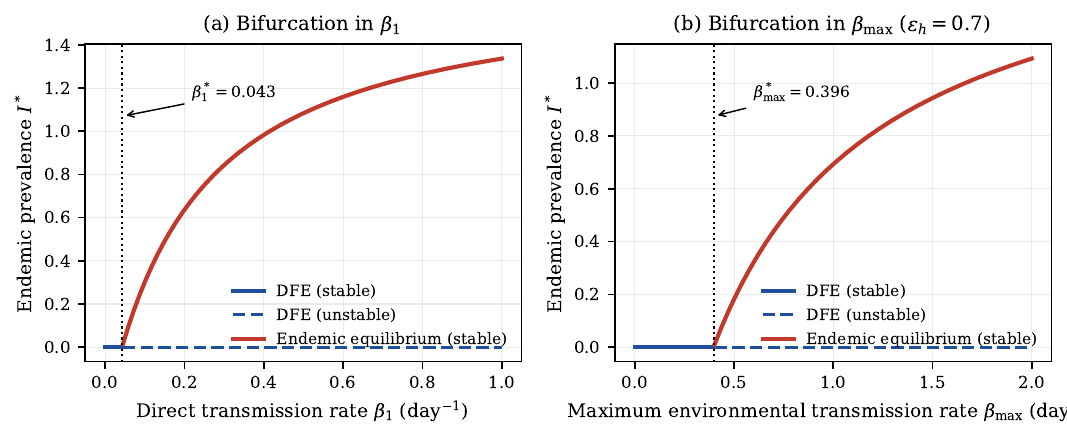}
\caption{Equilibrium continuation for \eqref{M1}, confirming \cref{thm:bifurcation}. Solid blue: stable disease-free branch; dashed blue: unstable disease-free branch; solid red: stable endemic branch, with stability determined from the eigenvalues of the full Jacobian at each point. \textbf{(a)} Continuation in $\beta_1$ at baseline (\cref{tab:baseline}); the exchange of stability occurs at $\beta_1^*=0.043$~day$^{-1}$ from \eqref{eq:beta1star}. \textbf{(b)} Continuation in $\beta_{\max}$ with $\epsilon_h=0.7$, so that $\Rzh=0.56<1$ and a threshold exists at $\beta_{\max}^*=0.396$~day$^{-1}$. In both panels, the endemic branch emerges continuously with no fold, and no branch of positive equilibria exists for $\Rz<1$. The biologically infeasible branch $I^*<0$ is omitted, as it lies outside $\Omega$; its counterpart in a backward bifurcation would be positive and unstable, producing the bistability shown here to be absent.}
\label{fig:bifurcation}
\end{figure}

\Cref{fig:bifurcation} presents numerical continuation of the equilibria of \eqref{M1}, obtained by solving \eqref{eq:EEscalar} and classifying stability from the eigenvalues of the full Jacobian. Panel (a) varies $\beta_1$ at baseline: the disease-free branch is stable for $\beta_1<\beta_1^*=0.043$ and loses stability at $\beta_1^*$, where a branch of stable endemic equilibria emerges continuously with zero jump and no fold, confirming \cref{thm:bifurcation}. Panel (b) repeats the analysis in $\beta_{\max}$. This requires $\Rzh<1$, since otherwise \eqref{eq:beta1star} has no analogue and $\Rz>1$ for every $\beta_{\max}$; we therefore fix $\epsilon_h=0.7$, giving $\Rzh=0.56$ and the threshold
\[
\beta_{\max}^*=\frac{k\sigma(\mu+\nu)A\,(1-\Rzh)}{\Lambda\theta(1-\epsilon_w)}=0.396~\text{day}^{-1},
\]
reproduced by the continuation to four significant figures. The same forward structure appears, showing that either pathway can drive the transition independently once the other is controlled.

Two features of the diagram deserve comment. First, the endemic prevalence is small in absolute terms: at baseline $I^*=1.08$ individuals, or $11$ per $100{,}000$, with $S^*=3{,}590$, $R^*=5{,}934$ and $W^*=32.8$ cells~mL$^{-1}$. This is not an artifact. With permanent immunity, the endemic state is sustained solely by the recruitment of new susceptibles at $\Lambda=0.457$ per day, so $I^*<\Lambda/A=1.69$ is necessarily true. The endemic equilibrium therefore describes a low-level residual transmission regime reached over decades, whereas the $100$-day dynamics of \cref{subsec:baseline,subsec:single} are a transient epidemic on a much faster timescale. Both are properties of the same system, and neither should be read as a prediction of the other.

Second, the negative branch of the transcritical bifurcation is omitted deliberately. For $\Rz<1$ the algebraic system \eqref{eq:EEsystem}--\eqref{eq:EEscalar} does possess a second root with $I^*<0$, which exchanges stability with the disease-free branch at $\Rz=1$ in the standard transcritical manner. That root lies outside $\mathbb{R}^4_{\ge0}$ and hence outside the feasible region $\Omega$, where,e by \cref{thm:invariance}, no trajectory of \eqref{M1} can reach it; it has no epidemiological interpretation, and plotting it would suggest a coexisting state that the model does not admit. Its existence is nonetheless what distinguishes forward from backward bifurcation: in a backward bifurcation, the second branch would be positive and unstable for $\Rz<1$, producing bistability and a fold. Its absence from the positive orthant, established analytically by $a<0$ in \eqref{eq:acoef} and by the monotonicity of $\mathcal{G}$, is precisely the property being demonstrated.

\subsection{Global sensitivity analysis}\label{subsec:sensitivity}

\subsubsection{Design}

Global sensitivity was assessed by Latin hypercube sampling with partial rank correlation coefficients \cite{Blower1994,Marino2008}. Three aspects of the design require justification.

\emph{Parameterisation.} As noted in \cref{subsec:params}, $\Lambda$ and $\mu$ are not independent inputs: sampling them independently over their marginal ranges implies community sizes $\Lambda/\mu$ spanning four orders of magnitude. At the same time, the initial condition is held fixed, so that the sampled systems are not comparable and $\Rze\propto\Lambda/\mu$ varies by a factor of $10^4$ for reasons unrelated to epidemiology. We therefore fix the community size at $N^\star=10^4$ and set $\Lambda=\mu N^\star$, sampling the remaining eleven parameters over the ranges of \cref{T1}. The sensitivity to population size is not lost, since by Proposition \ref{prop:elasticity} $\Upsilon_\Lambda=p_e=-\Upsilon_k$: population size and half-saturation concentration have equal and opposite elasticities, so the coefficient reported for $k$ also quantifies, with reversed sign, the effect of the size of the population served by the reservoir.

\emph{Scale.} Parameters spanning orders of magnitude ($\mu$, $k$, $\theta$) were sampled log-uniformly and the remainder uniformly. Uniform sampling of $k$ over $[10^5,10^7]$ would place $90\%$ of the mass above $10^6$, effectively removing the low-dose regime from the design.

\emph{Output.} The output is the attack ratio, that is, cumulative infections over $100$ days divided by $N^\star$, which is scale-free and does not saturate as readily as peak prevalence, whose lower bound is the seeded prevalence. We used $N=2{,}000$ samples; significance was assessed using a two-sided $t$ test on the coefficients with $N-d-1=1{,}988$ degrees of freedom, and $95\%$ confidence intervals were estimated from $300$ bootstrap resamples.

\subsubsection{Results}

\begin{figure}[!htbp]
\centering
\includegraphics[width=0.65\textwidth]{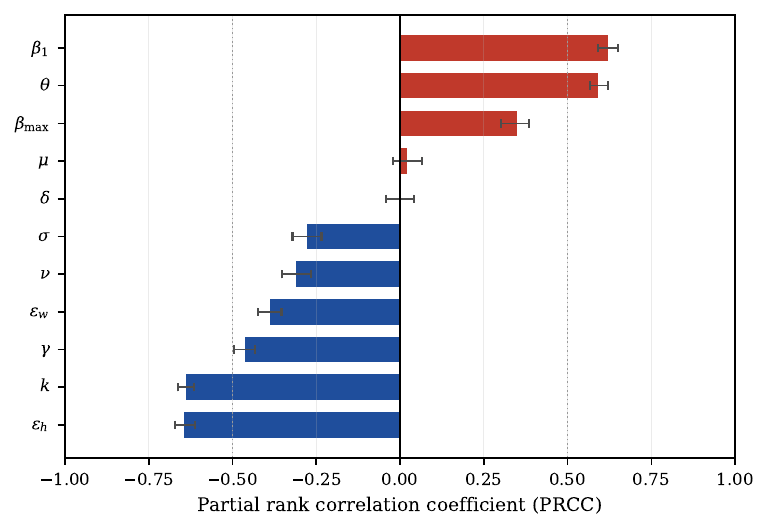}
\caption{Partial rank correlation coefficients for the $100$-day attack ratio, from $N=2{,}000$ Latin hypercube samples over the ranges in \cref{T1} with $\Lambda=\mu N^\star$, $N^\star=10^4$; $\mu$, $k$ and $\theta$ sampled log-uniformly. Bars are coloured by sign and error bars are $95\%$ bootstrap confidence intervals over $300$ resamples. All coefficients except those of $\mu$ ($p=0.3$) and $\delta$ ($p=0.9$) are significant at $p<10^{-30}$. Every sign agrees with the closed-form elasticities of Proposition \ref{prop:elasticity}; in particular, the half-saturation constant $k$ has a strongly negative coefficient, as $\Upsilon_k=-p_e<0$ requires.}
\label{sensitivity}
\end{figure}

\Cref{sensitivity} reports the coefficients. Human sanitation efficacy ($-0.65$), the half-saturation constant ($-0.64$), the direct transmission rate ($+0.62$) and the shedding rate ($+0.59$) are the strongest influences, followed by the recovery rate ($-0.46$), environmental sanitation efficacy ($-0.39$), maximal environmental transmission ($+0.35$), vaccination rate ($-0.31$) and pathogen removal rate ($-0.28$). Natural mortality ($+0.02$, $p=0.3$) and disease-induced mortality ($-0.00$, $p=0.9$) are indistinguishable from zero.

Every sign agrees with the analytic elasticities of Proposition \ref{prop:elasticity}, which provides an independent check on the implementation. In particular, the half-saturation constant has a strongly \emph{negative} coefficient. This is required: $k$ enters the model only through $\beta_2(W)=\beta_{\max}W/(k+W)$, and increasing it raises the concentration needed for a given force of infection, so $\Rze\propto1/k$ and $\Upsilon_k=-p_e<0$ exactly. There is no mechanism in \eqref{M1} by which a larger half-saturation constant could increase burden, so a positive coefficient for $k$ could only arise from an implementation error. We stress the point because the analytic elasticities make such errors detectable: whenever a computed index contradicts \eqref{eq:elas1}--\eqref{eq:elas3}, it is the computation rather than the intuition that is at fault. We recommend this check as routine practice whenever a reproduction number can be expressed in closed form. The corresponding interpretation is also mechanistically informative: the pair $(k,\theta)$ has nearly equal and opposite coefficients because burden depends on the reservoir essentially through the ratio $\theta/k$, the shedding per case measured in units of the infectious dose.

Two further points concern interpretation. First, the vaccination rate is influential but not dominant, in contrast to what is obtained when $\nu$ is sampled over $[0,0.1]$~day$^{-1}$; a partial rank correlation coefficient measures the contribution of a parameter to output variance \emph{given the assumed range}, and inflating the range of a control inflates its apparent importance. Our range is bounded by documented oral cholera vaccine delivery capacity \cite{azman2015impact}. Second, the negligible coefficient of $\mu$ confirms the timescale argument of Remark \ref{rem:timescale}: although $\Upsilon_\mu$ is not small, natural mortality cannot influence a $100$-day outbreak when $\mu^{-1}\approx60$ years. The two analyses are complementary rather than contradictory: elasticity ranks mechanisms at the threshold, partial rank correlation ranks uncertainty contributions to burden over a finite horizon, and disagreement between them locates a timescale separation rather than an error.

\subsection{Elimination boundaries in the sanitation plane}\label{subsec:R0control}

\begin{figure}[!htbp]
\centering
\includegraphics[width=0.55\textwidth]{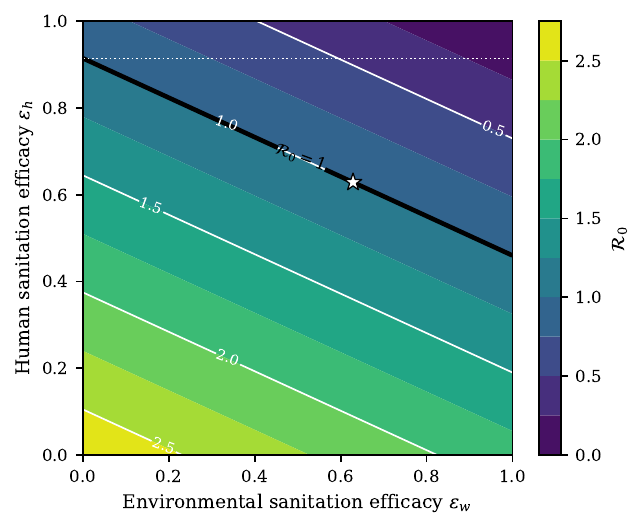}
\caption{Control reproduction number $\Rz=\Rzh+\Rze$ over the sanitation plane at $\nu=0$ and otherwise baseline values (\cref{tab:baseline}). Contours are straight because $\Rz$ is affine in $(\epsilon_h,\epsilon_w)$, with slope $-p_h/p_e=-2.2$ set by the pathway weights. The heavy contour is $\Rz=1$; the region above and to the right of it is the elimination regime. The dotted line marks $\epsilon_h=0.914$, the threshold for human sanitation used alone; environmental sanitation alone never reaches the boundary because $\Rzh=1.85>1$. The star marks equal joint efficacies $\epsilon_h=\epsilon_w=0.629$, the least demanding symmetric strategy achieving $\Rz=1$.}
\label{fig:R0control}
\end{figure}

\Cref{fig:R0control} shows $\Rz$ over the $(\epsilon_w,\epsilon_h)$ plane at $\nu=0$, computed from \eqref{eq:R0}. The contours are straight lines because $\Rz$ is affine in the efficacies, with slope $-p_h/p_e=-2.2$: the exchange rate between the two controls is fixed by the pathway weights, and one point of environmental efficacy substitutes for $0.45$ points of human efficacy at this baseline. The elimination boundary is
\[
p_h\epsilon_h+p_e\epsilon_w>1-\frac{1}{\Rz}=0.629 .
\]
Human sanitation alone requires $\epsilon_h>0.914$; environmental sanitation alone cannot cross the boundary at any efficacy, since $\Rzh=1.85>1$; but equal joint efficacies of $\epsilon_h=\epsilon_w=0.629$ suffice. The practical content is that integration converts an unattainable target into an attainable one, and Proposition \ref{prop:synergy} identifies exactly why: the boundary is a straight line rather than the convex curve that multiplicative independence would produce.

\section{Discussion}\label{sec5}

This study shows how the additive structure of cholera transmission can be used to connect epidemic thresholds, parameter sensitivity, and intervention effects within a parsimonious SIWR framework. The decomposition
\[
\Rz=\Rzh+\Rze
\]
is particularly informative because it separates direct and environmental transmission and determines their relative contributions to the sensitivity of the control reproduction number. The pathway weights $p_h$ and $p_e$ therefore provide a natural basis for comparing sanitation strategies across epidemiological settings. In our baseline parameterization, $p_h=0.69$, indicating greater sensitivity to human sanitation than to environmental sanitation at equal efficacy. This conclusion is setting dependent: a larger environmental contribution would shift the relative benefit toward environmental control. Thus, intervention priorities should reflect the underlying composition of cholera transmission rather than be interpreted as universal properties of the disease \cite{tien2010multiple,mukandavire2011estimating}.

The stability analysis also clarifies the relationship between vaccination and epidemic thresholds. \Cref{thm:global_DFE} establishes global stability under the sharpened condition $\Rtil\le1$, while Proposition \ref{prop:sharp} shows that the remaining gap with $\Rz<1$ cannot be removed using the present class of linear Lyapunov functions. This distinction reflects the transient difference between the initial susceptible population and the vaccinated disease-free equilibrium. Relatedly, the reproduction number evaluated at the vaccinated equilibrium reflects a long-term demographic state and should not be uncritically interpreted as the invasion threshold during a rapidly implemented vaccination campaign. The effective reproduction number $\Reff(t)$ provides a more appropriate complementary measure for such transient dynamics.

The analysis further provides a quantitative interpretation of combined sanitation. Under the Bliss independence criterion, the interaction between human and environmental sanitation contains a pathway component proportional to $p_hp_e\epsilon_h\epsilon_w$, together with an additional contribution arising from the nonlinear relationship between epidemic burden and the reproduction number. Consequently, synergy depends on both pathway composition and intervention intensity and need not increase monotonically with control effort. This explains why moderate combined interventions can produce substantial excess reductions in disease burden, whereas the measured excess benefit becomes small once transmission is strongly suppressed.

The model does not explicitly distinguish asymptomatic infection, although asymptomatic cases are common in cholera \cite{king2008inapparent,nelson2009transmission}. Their contribution is represented implicitly through effective transmission and shedding parameters, but this simplification may underestimate total infection burden and environmental transmission when calibration relies primarily on reported symptomatic cases. An explicit asymptomatic compartment, preferably informed by serological data, would therefore be a useful extension. More generally, the present results complement earlier SIWR studies by linking pathway-specific reproduction numbers to intervention sensitivity and interaction \cite{codeco2001endemic,tien2010multiple,mukandavire2011estimating}, while clarifying the interpretation of vaccination over epidemic and demographic timescales \cite{azman2015impact}.

Several limitations define directions for further work. The model assumes homogeneous mixing, constant intervention efficacy, permanent immunity over the epidemic horizon, and an autonomous environmental reservoir: spatial heterogeneity, time-dependent controls, waning immunity, and climate-driven variation in \emph{V.~cholerae} abundance could alter transient or long-term dynamics and warrant separate analysis \cite{alsammani2025impact,azman2013incubation,abdulrashid2019stability,Koelle2005,Pascual2002}. In particular, the forward-bifurcation result established here should not be extrapolated to models with additional mechanisms, such as waning immunity, that may change the equilibrium structure.

\section{Conclusion}\label{sec6}

This study developed a unified analytical framework for evaluating cholera transmission and control within a parsimonious SIWR model that incorporates human and environmental sanitation, as well as vaccination. The decomposition $\Rz=\Rzh+\Rze$ links the relative contributions of direct and environmental transmission to parameter sensitivity and intervention effectiveness. We established local and global stability conditions for the disease-free equilibrium. We proved that the bifurcation at $\Rz=1$ is forward for all biologically admissible parameter values, excluding backward bifurcation within the proposed model. The analysis further provides an explicit characterization of the interaction between sanitation measures under the Bliss independence criterion, while sensitivity analysis identifies the principal drivers of epidemic burden.

These results emphasize that effective cholera control depends on both the composition of transmission pathways and the coordinated use of interventions. In particular, combining controls targeting direct and environmental transmission can yield greater reductions in epidemic burden than either intervention alone. More broadly, the framework connects epidemic thresholds, pathway-specific sensitivity, and intervention interactions in a mathematically transparent setting, providing a basis for extensions that incorporate asymptomatic infection, spatial heterogeneity, seasonal forcing, and other mechanisms relevant to environmentally transmitted diseases.

\section*{Data Availability}
No new datasets were generated during this study. Numerical simulations were performed using published parameter values cited in the manuscript. Simulation code is available from the corresponding author upon reasonable request.

\section*{Conflict of Interest}
The authors declare that they have no conflicts of interest.

\end{document}